\documentclass[11pt,reqno]{amsproc}

\title[]{Critical SQG in bounded domains}

\author{Peter Constantin}
\address{Department of Mathematics, Princeton University, Princeton, NJ 08544}
\email{const@math.princeton.edu}

\author{Mihaela Ignatova}
\address{Department of Mathematics, Princeton University, Princeton, NJ 08544}
\email{ignatova@math.princeton.edu}

\usepackage[margin=1in]{geometry}
\usepackage{amsmath, amsthm, amssymb}
\usepackage{times}
\usepackage{color}
\usepackage{hyperref}

\newcommand{\pa}{\partial}
\newcommand{\la}{\label}
\newcommand{\fr}{\frac}
\newcommand{\na}{\nabla}
\newcommand{\be}{\begin{equation}}
\newcommand{\ee}{\end{equation}}
\newcommand{\ba}{\begin{array}{l}}
\newcommand{\ea}{\end{array}}
\newcommand{\Rr}{{\mathbb R}}

\newcommand{\beg}{\begin}

\newcommand{\D}{\Delta}
\renewcommand{\l}{\Lambda_D}
\date{today}
\begin{document}
\begin{abstract}
We consider the critical dissipative SQG equation in bounded domains, with the square root of the Dirichlet Laplacian dissipation. We prove global a priori interior $C^{\alpha}$ and Lipschitz bounds for large data. 
\end{abstract}


\noindent

\maketitle

\section{Introduction}

The Surface Quasigeostrophic equation (SQG) of geophysical origin (\cite{held}) was proposed as a two dimensional model for the study of inviscid incompressible formation of singularities (\cite{c}, \cite{cmt}). While the global regularity of all solutions of SQG whose initial data are smooth is still unknown, the original blow-up scenario of \cite{cmt} has been ruled out analytically (\cite{cord}) and numerically (\cite{cnum}), and nontrivial examples of global smooth solutions have been constructed (\cite{cascor}). Solutions of SQG and related equations without dissipation and with non-smooth (piece-wise constant) initial data give rise to interface dynamics (\cite{fr}, \cite{castro}) with potential finite time blow up (\cite{cornum}).

The addition of fractional Laplacian dissipation produces globally regular solutions if the power of the Laplacian is larger or equal than one half. When the linear dissipative operator is precisely the square root of the Laplacian, the equation is commonly referred to as the ``critical dissipative SQG'', or ``critical SQG''. This active scalar equation (\cite{c})  has been the object of intensive study in  the past decade. The solutions are transported by divergence-free velocities they create, and are smoothed out and decay due to  nonlocal diffusion. Transport and diffusion do not add size to a solution: the solution remains bounded, if it starts so (\cite{res}). The space $L^{\infty}(\Rr^2)$ is not a natural phase space for the nonlinear evolution: the nonlinearity involves Riesz transforms and these are not well behaved in $L^{\infty}$. Unfortunately, for the purposes of studies of global in time behavior of solutions, $L^{\infty}$ is unavoidable: it quantifies the most important information freely available. The equation is quasilinear and $L^{\infty}$--critical, and there is no `` wiggle room'', nor a known better (smaller) space which is invariant for the evolution. One must work in order to obtain better information. A pleasant aspect of criticality is that solutions with small initial $L^{\infty}$ norm are smooth (\cite{ccw}). The global regularity of large solutions was obtained independently in \cite{caf} and \cite{knv} by very different methods: using harmonic extension and the De Georgi methodology  of zooming in, and passing from $L^2$ to $L^{\infty}$ and from $L^{\infty}$ to $C^{\alpha}$ in \cite{caf}, and constructing a family of time-invariant moduli of continuity in \cite{knv}. Several subsequent proofs were obtained (please see \cite{cvt} and references therein). All the proofs are dimension-independent, but are in either $\Rr^d$ or on the torus ${\mathbb {T}}^d$. The proofs of \cite{cv1} and \cite{cvt} were based on an extension of the C\'{o}rdoba-C\'{o}rdoba inequality (\cite{cc}).
This inequality states that
\be
\Phi'(f)\Lambda f - \Lambda\Phi(f)\ge 0
\la{corcorone}
\ee
pointwise. Here $\Lambda = \sqrt{-\Delta}$ is the square root of the Laplacian in the whole space $\Rr^d$, $\Phi$ is a real valued convex function of one variable, normalized so that $\Phi(0) = 0$ and $f$ is a smooth function. The fractional Laplacian in the whole space has a (very) singular integral representation, and this can be used to obtain (\ref{corcorone}). In \cite{cv1} specific nonlinear maximum principle lower bounds were obtained and used to prove the global regularity. A typical example is
\be
D(f) =f\Lambda f- \fr{1}{2}\Lambda\left({f^2}\right) \ge c\left(\|\theta\|_{L^{\infty}}\right)^{-1} {f^3}
\la{corcorv}
\ee
pointwise, for $f=\pa_i\theta$ a component of the gradient of a bounded function $\theta$. This is a useful cubic lower bound for a quadratic expression, when  $\|\theta\|_{L^{\infty}}\le\|\theta_0\|_{L^{\infty}}$ is known to be bounded above. The critical SQG equation in $\Rr^2$ is 
\be
\pa_t \theta + u\cdot\na\theta + \Lambda\theta = 0
\la{critsqg}
\ee
where
\be
u = \na^{\perp}\Lambda^{-1}\theta  = R^{\perp}{\theta}
\la{ucr}
\ee
and $\na^{\perp} = (-\pa_2, \pa_1)$ is the gradient rotated by $\fr{\pi}{2}$.
Because of the conservative nature of transport and the good dissipative properties of $\Lambda$ following from (\ref{corcorone}), all $L^p$ norms of $\theta$ are nonincreasing in time. Moreover, because of properties of Riesz transforms, $u$ is essentially of the same order of magnitude as $\theta$. Differentiating the equation we obtain the stretching equation
\be
\left(\pa_t + u\cdot\na + \Lambda\right)\na^{\perp}\theta = (\na u)\na^{\perp}\theta.
\la{stretch}
\ee
(In the absence of $\Lambda$ this is the same as the stretching equation for three dimensional vorticity in incompressible Euler equations, one of the main reasons SQG was considered in 
\cite{c}, \cite{cmt} in the first place.) Taking the scalar product with $\na^{\perp}\theta$ we obtain
\be
\fr{1}{2}(\pa_t + u\cdot\na + \Lambda)q^2 + D(q) = Q
\la{qf}
\ee
for $q^2 = |\na^{\perp}\theta|^2$, with
\[
Q = (\na u)\na^{\perp}\theta\cdot\na^{\perp}\theta \le |\na u| q^2.
\]
The operator $\pa_t + u\cdot\na + \Lambda$ is an operator of advection and fractional diffusion: it does not add size. Using the pointwise  bound (\ref{corcorv}) we already see that the dissipative lower bound is potentially capable of dominating the cubic term $Q$, but there are two obstacles. The first obstacle is that constants matter: the two expressions are cubic, but the useful dissipative cubic lower bound $D(q)\ge K |q|^{3}$ has perhaps too small a prefactor $K$ if the $L^{\infty}$ norm of $\theta_0$ is too large.  The second obstacle is that although
\[
\na u = R^{\perp}(\na\theta)
\]
has the same size as $\na^{\perp}\theta$ (modulo constants) in all $L^p$ spaces $1<p<\infty$, it fails to be bounded in $L^{\infty}$ by the $L^{\infty}$ norm of $\na^{\perp}\theta$. In order to overcome these obstacles, in \cite{cv1} and \cite{cvt}, instead of estimating directly gradients, the proof proceeds by estimating finite differences, with the aim of obtaining bounds for $C^{\alpha}$ norms first. In fact, in critical SQG, once
the solution is bounded in any $C^{\alpha}$ with $\alpha>0$, it follows that it is $C^{\infty}$. More generally, if the equation has a dissipation of order $s$, i.e., $\Lambda$ is replaced by $\Lambda^s$ with $0<s\le 1$ then if $\theta$ is bounded in $C^{\alpha}$ with $\alpha>1-s$, then the solution is smooth (\cite{cw}). (This condition is sharp, if one considers  general linear advection diffusion equations, (\cite{svz}). In \cite{cvt} the smallness of $\alpha$ is used to show that the term corresponding to $Q$ in the finite difference version of the argument is dominated by the term corresponding to $D(q)$.

In this paper we consider the critical SQG equation in bounded domains.  We take a bounded open domain $\Omega\subset \Rr^d$ with smooth (at least $C^{2,\alpha}$) boundary and denote by $\D$ the Laplacian operator with homogeneous Dirichlet boundary conditions and by $\l$ its square root defined in terms of eigenfunction expansions. Because no explicit kernel for the fractional Laplacian is available in general, our approach, initiated in \cite{ci} is based on bounds on the heat kernel.

The critical SQG equation is 
\be
\pa_t \theta + u\cdot\na\theta + \l\theta =0
\la{sqg}
\ee
with
\be
u = \na^{\perp}\l^{-1}\theta = R_D^{\perp}\theta
\la{u}
\ee
and smooth initial data.  We obtain global regularity results, in the spirit of the ones in the whole space. There are quite significant differences between the two cases. First of all, the fact that no explicit formulas are available for kernels requires a new approach; this yields as a byproduct new proofs even in the whole space. The main difference and additional difficulty  in the bounded domain case is due to the lack of translation invariance.
The fractional Laplacian is not translation invariant, and from the very start, differentiating the equation (or taking finite differences) requires understanding the respective commutators. For the same reason, the Riesz transforms $R_D$ are not spectral operators, i.e., they do not commute with functions of the Laplacian, and so velocity bounds need a different treatment.  In \cite{ci} we proved using the heat kernel approach the existence of global weak solutions of (\ref{sqg}) in $L^2(\Omega)$. A proof of local existence of smooth solutions is provided in the present paper in $d=2$. The local existence is obtained in Sobolev spaces based on $L^2$ and uses Sobolev embeddings. Because of this, the proof is dimension dependent. A proof in higher dimensions is also possible but we do not pursue this here. We note that for regular enough solutions (e.g. $\theta\in H_0^1(\Omega)$) the normal component of the velocity 
vanishes at the boundary $\left(R_D^{\perp}\theta\cdot N\right)_{\left |\right. \pa\Omega}=0$ because the stream function $\psi = \l^{-1}\theta$ vanishes at the boundary and its gradient is normal to the boundary. Let us remark here that even in the case of a half-space and $\theta\in C_0^{\infty}(\Omega)$, the tangential component of the velocity need not vanish: there is tangential slip.

In order to state our main results, let 
\be
d(x) = dist(x,\pa\Omega)
\la{dx}
\ee
denote the distance from $x$ to the boundary of $\Omega$. We introduce the $C^{\alpha}(\Omega)$ space for interior estimates:
\beg{defi}\la{calpha}
Let $\Omega$ be a bounded domain and let  $0<\alpha<1$ be fixed. We say that $\theta\in C^{\alpha}(\Omega)$ if $\theta\in L^{\infty}(\Omega)$ and
\be
[f]_{\alpha} = \sup_{x\in\Omega}(d(x))^{\alpha}\left(\sup_{h\neq 0,|h|<d(x)}\fr{|f(x+h)-f(x)|}{|h|^{\alpha}}\right) <\infty.
\la{semi}
\ee
The norm in $C^{\alpha}(\Omega)$ is 
\be
\|f\|_{C^{\alpha}} = \|f\|_{L^{\infty}(\Omega)} + [f]_{\alpha}.
\la{norm}
\ee
\end{defi}
Our main results are the following:
\beg{thm}\la{alphaint} Let $\theta(x,t)$ be a smooth solution of (\ref{sqg}) on a time interval $[0, T)$, with $T\le \infty$, with initial data $\theta(x,0)= \theta_0(x)$. Then the solution is uniformly bounded,
\be
\sup_{0\le t< T}\|\theta(t)\|_{L^{\infty}(\Omega)}\le \|\theta_0\|_{L^{\infty}(\Omega)}.
\la{linftyb}
\ee
There exists $\alpha$ depending only on $\|\theta_0\|_{L^{\infty}(\Omega)}$ and
$\Omega$, and a constant $\Gamma$ depending only on the domain $\Omega$ (and in particular, independent of $T$) such that
\be
\sup_{0\le t<T}\|\theta(t)\|_{C^{\alpha}(\Omega)} \le \Gamma\|\theta_0\|_{C^{\alpha}(\Omega)} 
\la{calphain}
\ee
holds.
\end{thm}
The second theorem is about global interior gradient bounds: 
\beg{thm}\la{gradint} Let $\theta(x,t)$ be a smooth solution of (\ref{sqg}) on a time interval $[0, T)$, with $T\le \infty$, with initial data $\theta(x,0)= \theta_0(x)$. There exists a constant $\Gamma_1$ depending only on $\Omega$ such that
\be
\sup_{x\in\Omega, 0\le t<T}d(x)|\na_x\theta(x,t)|\le \Gamma_1\left[\sup_{x\in\Omega}d(x)|\na_x\theta_0(x)| + \left(1+\|\theta_0\|_{L^{\infty}(\Omega)}\right)^4\right]
\la{gradintb}
\ee
holds.
\end{thm}
\beg{rem} Higher interior regularity can be proved also. In fact, once global interior $C^{\alpha}$ bounds are obtained for any $\alpha>0$, the interior regularity problem becomes subcritical, meaning that ``there is room to spare''. This is already the case for Theorem \ref{gradint} and justifies thinking that the equation is $L^{\infty}$ interior-critical. However, we were not able to obtain global uniform $C^{\alpha}(\bar{\Omega})$ bounds. Moreover, we do not know the implication 
$C^{\alpha}(\bar{\Omega}) \Rightarrow C^{\infty}(\bar{\Omega})$ uniformly, and thus the equation is not $L^{\infty}$ critical up to the boundary. This is due to the fact that the commutator between normal derivatives and the fractional Dirichlet Laplacian is not controlled uniformly up to the boundary. The example of half-space is instructive because explicit kernels and calculations are available. In this example odd reflection across the boundary permits the construction of global smooth solutions, if the initial data are smooth and compactly supported away from the boundary. The support of the solution remains compact and  cannot reach the boundary in finite time, but the gradient of the solution might grow in time at an exponential rate.
\end{rem} 

The proofs of our main results use the following elements. First, the inequality (\ref{corcorone}) which has been proved in (\cite{ci}) for the Dirichlet $\l$ is shown to have a lower bound 
\be
D(f)(x) = \left(f\l f - \fr{1}{2}\l\left({f^2}\right)\right)(x) \ge c\fr{f^2(x)}{d(x)}
\la{dfdxb}
\ee
with $c>0$ depending only on $\Omega$. Note that in $\Rr^d$,  $d(x)=\infty$, which is consistent with (\ref{corcorone}). This lower bound (valid for general $\Phi$ convex, with $c$ independent of $\Phi$, see (\ref{cor})) provides a strong damping boundary repulsive term, which is essential to overcome boundary effects coming from the lack of translation invariance.

The second element of proofs consists of  nonlinear lower bounds in the spirit of (\cite{cv1}).  A version for derivatives in bounded domains, proved in (\cite{ci}) is modified for finite differences. In order to make sense of finite differences near the boundary in a manner suitable for transport, we introduce a family of good cutoff functions depending on a scale $\ell$ in Lemma \ref{goodcutoff}. The finite difference nonlinear lower bound is
\be
D(f)(x)\ge c\left(|h|\|\theta\|_{L^{\infty}(\Omega)}\right)^{-1}|f(x)|^3+ c\fr{|f(x)|^2}{d(x)}
\la{dfcube}
\ee 
when $f=\chi\delta_h\theta$ is large (see (\ref{nlb})), where $\chi$ belongs to the family of good cutoff functions. 

Once global interior $C^{\alpha}(\Omega)$ bounds are obtained, in order to obtain global interior bounds for the gradient, we use a different nonlinear lower bound,
\be
D(f) = \left(f\l f -\fr{1}{2}(\l f^2)\right)(x) \ge c \fr{|f(x)|^{3+\fr{\alpha}{1-\alpha}}}{\|\theta\|_{C^{\alpha}(\Omega)}^{\fr{1}{1-\alpha}}}(d(x))^{\fr{\alpha}{1-\alpha}} + c\fr{f^2(x)}{d(x)}
\la{nnlbd}
\ee
for large $f=\chi\na\theta$ (see (\ref{nlbd})). This is a super-cubic bound, and makes the gradient equation look subcritical. Similar bounds were obtained in the whole space in (\cite{cv1}). Proving the bounds (\ref{dfcube}) and (\ref{nnlbd}) requires a different approach and new ideas because of the absence of explicit formulas and lack of translation invariance.

The third element of proofs are bounds for $R_D^{\perp}\theta$ based  only on global apriori information on $\|\theta\|_{L^{\infty}}$ and the nonlinear lower bounds on $D(f)$ for appropriate $f$.  Such an approach was initiated in (\cite{cv1}) and (\cite{cvt}). In the bounded domain case, again, the method of proof is different because the kernels are not explicit, and reference is made to the heat kernels. The boundaries introduce additional error terms. The bound for finite differences is
\be
|\delta_h R^{\perp}_D\theta(x)| \le  C\left(\sqrt{\rho D(f)(x)} + \|\theta\|_{L^{\infty}}\left(\fr{|h|}{d(x)}+ \fr{|h|}{\rho}\right) + |\delta_h\theta(x)|\right)
\la{dhrtb}
\ee
for  $\rho\le cd(x)$, with $f=\chi\delta_h \theta$ and  with $C$ a constant depending on $\Omega$ (see \ref{dhub}). The bound for gradient is
\be
|\na R^{\perp}_D\theta(x)| \le C\left(\sqrt{\rho D(f)(x)} + \|\theta\|_{L^{\infty}(\Omega)}\left(\fr{1}{d(x)} + \fr{1}{\rho}\right) + |\na\theta(x)|\right)
\la{nartxb}
\ee
for $\rho\le cd(x)$ with $f=\chi\na\theta$ with a constant $C$ depending on $\Omega$ (see (\ref{nauxb})). These are remarkable pointwise bounds (clearly not valid for the case of the Laplacian even in the whole space, where $D(f)(x) = |\na f(x)|^2$).

The fourth element of the proof are bounds for commutators. These bounds 
\be
\left |\left[\chi\delta_h, \l\right]\theta(x)\right| \le C\fr{|h|}{d(x)^2}\|\theta\|_{L^{\infty}(\Omega)},
\la{chidellcomm}
\ee
for $\ell\le d(x)$, (see (\ref{commhb})), and
\be
\left|\left[\chi\na, \l\right]\theta(x)\right|  \le \fr{C}{d(x)^2}\|\theta\|_{L^{\infty}(\Omega)},
\la{nachicommd}
\ee
for $\ell\le d(x)$, (see (\ref{nachi})), reflect the difficulties due to the boundaries. They are remarkable though in that the only price to pay for a second order commutator in $L^{\infty}$ is $d(x)^{-2}$. Note that in the whole space this commutator vanishes ($\chi=1$). This nontrivial situation in bounded domains is due to cancellations and bounds on the heat kernel representing translation invariance effects away from boundaries (see (\ref{cancel1}, \ref{cancel2})). Although the heat kernel in bounded domains has been extensively studied, and the proofs of (\ref{cancel1}) and (\ref{cancel2}) are elementary, we have included them in the paper because we have not found them readily available in the literature and for the sake of completeness.

The paper is organized as follows: after preliminary background, we
prove the nonlinear lower bounds. We have separate sections for bounds for the Riesz transforms and the commutators. The proof of the main results are then provided, using nonlinear maximum principles. We give some of the explicit calculations in the example of a half-space and conclude the paper by proving the translation invariance bounds for the heat kernel (\ref{cancel1}), (\ref{cancel2}), and a local well-posedness  result in two appendices.

\section{Preliminaries}
The $L^2(\Omega)$ - normalized eigenfunctions of $-\D$ are denoted $w_j$, and its eigenvalues counted with their multiplicities are denoted $\lambda_j$: 
\be
-\D w_j = \lambda_j w_j.
\la{ef}
\ee
It is well known that $0<\lambda_1\le...\le \lambda_j\to \infty$  and that $-\D$ is a positive selfadjoint operator in $L^2(\Omega)$ with domain ${\mathcal{D}}\left(-\D\right) = H^2(\Omega)\cap H_0^1(\Omega)$.
The ground state $w_1$ is positive and
\be
c_0d(x) \le w_1(x)\le C_0d(x)
\la{phione}
\ee
holds for all $x\in\Omega$, where $c_0, \, C_0$ are positive constants depending on $\Omega$. Functional calculus can be defined using the eigenfunction expansion. In particular
\be
\left(-\D\right)^{\beta}f = \sum_{j=1}^{\infty}\lambda_j^{\beta} f_j w_j
\la{funct}
\ee
with 
\[
f_j =\int_{\Omega}f(y)w_j(y)dy
\]
for $f\in{\mathcal{D}}\left(\left (-\D\right)^{\beta}\right) = \{f\left |\right. \; (\lambda_j^{\beta}f_j)\in \ell^2(\mathbb N)\}$.
We will denote by
\be
\l^s = \left(-\D\right)^{\fr{s}{2}}, 
\la{lambdas}
\ee
the fractional powers of the Dirichlet Laplacian, with $0\le s \le 2$
and with $\|f\|_{s,D}$ the norm in ${\mathcal{D}}\left (\l^s\right)$:
\be
\|f\|_{s,D}^2 = \sum_{j=1}^{\infty}\lambda_j^{s}f_j^2.
\la{norms}
\ee
It is well-known and easy to show that
\[
{\mathcal{D}}\left( \l \right) = H_0^1(\Omega).
\]
Indeed, for $f\in{\mathcal{D}}\left (-\D\right)$ we have
\be
\|\na f\|^2_{L^2(\Omega)} = \int_{\Omega}f\left(-\D\right)fdx = \|\l f\|_{L^2(\Omega)}^2 = \|f\|^2_{1,D}. 
\la{kat}
\ee
We recall that the Poincar\'{e} inequality implies that the Dirichlet integral on the left-hand side above is equivalent to the norm in $H_0^1(\Omega)$
and therefore the identity map from the dense subset ${\mathcal{D}}\left(-\D\right)$ of $H_0^1(\Omega)$ to ${\mathcal D}\left(\l\right)$ is an isometry, and thus $H_0^1(\Omega)\subset {\mathcal{D}}\left(\l\right)$. But ${\mathcal{D}}\left(-\D\right)$ is dense in  ${\mathcal D}\left(\l\right)$ as well, because finite linear combinations of eigenfunctions are dense in  ${\mathcal D}\left(\l\right)$. Thus the opposite inclusion is also true, by the same isometry argument.  \\
Note that in view of the identity
\be
\lambda^{\fr{s}{2}} = c_{s}\int_0^{\infty}(1-e^{-t\lambda})t^{-1-\fr{s}{2}}dt,
\la{lambdalpha}
\ee
with 
\[
1 = c_{s} \int_0^{\infty}(1-e^{-\tau})\tau^{-1-\fr{s}{2}}d\tau,
\]
valid for $0\le s <2$, we have the representation
\be
\left(\left(\l\right)^{s}f\right)(x) = c_{s}\int_0^{\infty}\left[f(x)-e^{t\D}f(x)\right]t^{-1-\fr{s}{2}}dt
\la{rep}
\ee
for $f\in{\mathcal{D}}\left(\left (-\l\right)^{s}\right)$.
We use precise upper and lower bounds for the kernel $H_D(t,x,y)$ of the heat operator,
\be
(e^{t\D}f)(x) = \int_{\Omega}H_D(t,x,y)f(y)dy .
\la{heat}
\ee
These are as follows (\cite{davies1},\cite{qszhang1},\cite{qszhang2}).
There exists a time $T>0$ depending on the domain $\Omega$ and constants
$c$, $C$, $k$, $K$, depending on $T$ and $\Omega $ such that
\be
\ba
c\min\left (\fr{w_1(x)}{|x-y|}, 1\right)\min\left (\fr{w_1(y)}{|x-y|}, 1\right)t^{-\fr{d}{2}}e^{-\fr{|x-y|^2}{kt}}\le \\H_D(t,x,y)\le C
\min\left (\fr{w_1(x)}{|x-y|}, 1\right)\min\left (\fr{w_1(y)}{|x-y|}, 1\right)t^{-\fr{d}{2}}e^{-\fr{|x-y|^2}{Kt}}
\ea
\la{hb}
\ee
holds for all $0\le t\le T$. Moreover
\be
\fr{\left |\na_x H_D(t,x,y)\right|}{H_D(t,x,y)}\le
C\left\{
\ba
\fr{1}{d(x)},\quad\quad \quad\quad {\mbox{if}}\; \sqrt{t}\ge d(x),\\
\fr{1}{\sqrt{t}}\left (1 + \fr{|x-y|}{\sqrt{t}}\right),\;{\mbox{if}}\; \sqrt{t}\le d(x)
\ea
\right.
\la{grbx}
\ee
holds for all $0\le t\le T$.
Note that, in view of
\be
H_D(t,x,y) = \sum_{j=1}^{\infty}e^{-t\lambda_j}w_j(x)w_j(y) ,
\la{hphi}
\ee
elliptic regularity estimates and Sobolev embedding which imply uniform absolute convergence of the series (if $\pa\Omega$ is smooth enough), we have that
\be
\pa_1^{\beta}H_D(t,y,x) = \pa_2^{\beta}H_D(t,x,y)
= \sum_{j=1}^{\infty}e^{-t\lambda_j}\pa_y^{\beta}w_j(y)w_j(x)
\la{dh}
\ee
for positive $t$, where we denoted by $\pa_1^{\beta}$ and $\pa_2^{\beta}$ derivatives with respect to the first spatial variables and the second spatial variables, respectively.

Therefore, the gradient bounds (\ref{grbx}) result in
\be
\fr{\left |\na_y H_D(t,x,y)\right|}{H_D(t,x,y)}\le
C\left\{
\ba
\fr{1}{d(y)},\quad\quad \quad\quad\quad {\mbox{if}}\; \sqrt{t}\ge d(y),\\
\fr{1}{\sqrt{t}}\left (1 + \fr{|x-y|}{\sqrt{t}}\right),\;{\mbox{if}}\; \sqrt{t}\le d(y).
\ea
\right.  
\la{grby}
\ee

We also use a bound
\be
\na_x\na_x H_D(x,y,t) \le Ct^{-1-\fr{d}{2}}e^{-\fr{|x-y|^2}{\tilde{K}t}}
\la{naxnaxb}
\ee
valid for $t\le cd(x)^2$ and $0<t\le T$, which follows from the upper bounds (\ref{hb}), (\ref{grbx}). 

Important additional bounds we need are 
\be
\int_{\Omega}\left |(\na_x +\na_y)H_D(x,y,t)\right|dy \le Ct^{-\fr{1}{2}}e^{-\fr{d(x)^2}{\tilde{K}t}}
\la{cancel1}
\ee
and
\be
\int_{\Omega}\left |\na_x(\na_x +\na_y)H_D(x,y,t)\right|dy \le Ct^{-1}e^{-\fr{d(x)^2}{\tilde{K}t}}
\la{cancel2}
\ee
valid for $t\le cd(x)^2$ and $0<t\le T$. These bounds reflect the fact that translation invariance is remembered in the solution of the heat equation with Dirichlet boundary data for short time, away from the boundary. We sketch the proofs of (\ref{naxnaxb}), (\ref{cancel1}) and (\ref{cancel2}) in the Appendix 1.

\section{Nonlinear Lower Bounds}
We prove bounds in the spirit of (\cite{cv1}). The proofs below are based on the method of (\cite{ci}), but they concern different objects (finite differences, properly localized) or different assumptions ($C^{\alpha}$). Nonlinear lower bounds are an essential ingredient in proofs of global regularity for drift-diffusion equations with nonlocal dissipation.

We start with a couple lemmas. In what follows we denote by $c$ and $C$ generic positive constants that depend on $\Omega$. When the logic demands it, we temporarily manipulate them and number them to show that the arguments are not circular. There is no attempt to optimize constants, and their numbering is local in the proof, meaning that, if for instance  $C_2$ appears in two proofs, it need not be the same constant. However, when emphasis is necessary we single out constants, but then we  avoid the letters $c,C$ with or without subscripts.

\beg{lemma}\la{Thet}
The solution of the heat equation with initial datum equal to 1 and zero boundary conditions,
\be
\Theta(x,t) = \int_{\Omega}H_D(x,y,t)dy
\la{Theta}
\ee
obeys $0\le \Theta(x,t)\le 1$, because  of the maximum principle. There exist constants $T, c, C$ depending only on $\Omega$ such that the following inequalities hold:
\be
\Theta(x,t)\ge c\min\left\{1, \left(\fr{d(x)}{\sqrt{t}}\right)^d\right\}
\la{thetalow}
\ee
for all $0\le t\le T$, and
\be
\Theta(x,t) \le C \fr{d(x)}{\sqrt{t}}
\la{thetaup}
\ee
for all $0\le t\le T$. Let $0<s<2$. There exists a constant $c$ depending on $\Omega$ and $s$ such that
\be
\int_0^{\infty}t^{-1-\fr{s}{2}}(1-\Theta(x,t))dt \ge c d(x)^{-s}
\la{lambdatheta}
\ee 
holds.

\end{lemma}
\beg{rem} $\l^{s} 1 $ is defined by duality by the left hand side of (\ref{lambdatheta})  and belongs to $H^{-1}(\Omega)$. 
\end{rem}
\noindent{\bf Proof.} Indeed,
\[
\Theta (x,t) = \int_{\Omega}H_D(t,x,y)dy\ge \int_{|x-y|\le \fr{d(x)}{2}} H_D(t,x,y)dy
\]
because $H_D$ is positive. Using the lower bound in (\ref{phione}) we have that
$|x-y|\le \fr{d(x)}{2}$ implies
\[
\fr{w_1(x)}{|x-y|}\ge 2c_0,\quad \fr{w_1(y)}{|x-y|}\ge c_0,
\]
and then, using the lower bound in (\ref{hb}) we obtain
\[
H_D(t,x,y) \ge 2cc_0^2t^{-\fr{d}{2}}e^{-\fr{|x-y|^2}{kt}}. 
\]
Integrating it follows that
\[
\Theta(x,t) \ge 2 cc_0^2\omega_{d-1}k^{\fr{d}{2}}\int_0^{\fr{d(x)}{2\sqrt{kt}}}\rho^{d-1}e^{-\rho^2}d\rho.
\]
If $\fr{d(x)}{2\sqrt{kt}}\ge 1$ then the integral is bounded below by
$\int_0^1\rho^{d-1}e^{-\rho^2}d\rho$. If $\fr{d(x)}{2\sqrt{kt}}\le 1$ then
$\rho\le 1$ implies that the exponential is bounded below by $e^{-1}$ and so
(\ref{thetalow}) holds.

Now (\ref{thetaup}) holds immediately from (\ref{phione}) and the upper bound 
in (\ref{hb}) because the integral
\[
\int_{\Rr^d}|\xi|^{-1}e^{-\fr{|\xi|^2}{K}}d\xi <\infty
\]
if $d\ge 2$.

Regarding (\ref{lambdatheta}) we use
\[
\int_0^{\infty}t^{-1-\fr{s}{2}}(1-\Theta(x,t))dt\ge \int_\tau^{T}t^{-1-\fr{s}{2}}(1-\Theta(x,t))dt
\]
and choose appropriately $\tau$. In view of (\ref{thetaup}), if 
\[
\fr{d(x)}{\sqrt{\tau}}\le \fr{1}{2C}
\]
then, when  $\tau\le t\le T$ we have
\[
1-\Theta(x,t)\ge \fr{1}{2},
\]
and therefore
\[
\int_{\tau}^{T} t^{-1-\fr{s}{2}}\left(1-\Theta(x,t)\right)dt \ge\fr{1}{s} \tau^{-\fr{s}{2}}\left ( 1- \left(\fr{\tau}{T}\right)^{\fr{s}{2}}\right)
\]
holds. The choice 
\[
\fr{d(x)}{\sqrt{\tau}} =\fr{1}{2C}
\]
implies (\ref{lambdatheta}) provided $2\tau \le T$ which is the same as $d(x)\le \fr{\sqrt{T}}{2C\sqrt{2}}$. On the other hand, $\Theta$ is exponentially small if $t$ is large enough, so the contribution to the integral in (\ref{lambdatheta}) is bounded below by a nonzero constant. This ends the proof of the lemma.

\beg{lemma}\la{grbnd}
Let $0\le \alpha<1$. There exists constant $C$ depending on $\Omega$ and 
$\alpha$ such that
\be
\int_{\Omega}|\na_y H_D(t,x,y)| |x-y|^{\alpha}dy \le C t^{-\fr{1-\alpha}{2}}
\la{grup}
\ee
holds for $0\le t\le T$.
\end{lemma}
Indeed, the upper bounds (\ref{hb}) and  (\ref{grby}) yield
\[
\ba
\int_{d(y)\ge \sqrt{t}}|\na_y H_D(t,x,y)||x-y|^{\alpha}dy \\\le C_2t^{-\fr{1}{2}}
\int_{\Rr^d}\left (1 + \fr{|x-y|}{\sqrt{t}}\right)t^{-\fr{d}{2}}e^{-\fr{|x-y|^2}{Kt}}|x-y|^{\alpha}dy\\
= C_3t^{-\fr{1-\alpha}{2}}
\ea
\]
and, in view of the upper bound in (\ref{phione}), $\fr{1}{d(y)}w_1(y)\le C_0$ and the upper bound in (\ref{hb}), we have
\[
\ba
\int_{d(y)\le \sqrt{t}}|\na_y H_D(t,x,y)||x-y|^{\alpha}dy \\\le C_4\int_{\Rr^d}\fr{1}{|x-y|}t^{-\fr{d}{2}}e^{-\fr{|x-y|^2}{Kt}}|x-y|^{\alpha}dy = C_5t^{-\fr{1-\alpha}{2}}.
\ea
\]
This proves (\ref{grup}). We introduce now a good family of cutoff functions $\chi$ depending on a length scale $\ell$.

\beg{lemma}\la{goodcutoff} Let $\Omega$ be a bounded domain with $C^2$ boundary. For $\ell>0$ small enough (depending on $\Omega$) there exist cutoff functions $\chi$ with the properties: $0\le \chi\le 1$, $\chi(y)=0$ if $d(y)\le \fr{\ell}{4}$, $\chi(y)= 1$ for $d(y)\ge \fr{\ell}{2}$, $|\na^k\chi|\le C\ell^{-k}$ with $C$ independent of $\ell$ and 
\be
\int_{\Omega}\fr{(1-\chi(y))}{|x-y|^{d+j}}dy \le C\fr{1}{d(x)^{j}}
\la{chij}
\ee
and
\be
\int_{\Omega}|\na\chi(y)|\fr{1}{|x-y|^{d-\alpha}}\le Cd(x)^{-(1-\alpha)}
\la{nachij}
\ee
hold for $j>-d$, $\alpha<d$ and $d(x)\ge \ell$. We will refer to such $\chi$ as a ``good cutoff''.
\end{lemma}
\noindent{\bf Proof.} There exists a length $\ell_0$ such that if $P$ is a point of the boundary $\pa\Omega$, and if $|P-y|\le 2\ell_0$, then $y\in\Omega$ if and only if (after a rotation and a translation) $y_d>F(y')$, where
$y'=(y_1, \dots, y_{d-1})$ and $F$ is a $C^2$ function with $F(0)=0$, $\na F(0)=0$, $|\na F|\le\fr{1}{10}$. We took thus without loss of generality coordinates such that $P= (0,0)$ and the normal to $\pa\Omega$ at $P$ is $(0,\dots, 0, 1)$. Now if $\ell<\ell_0$ and $d(x)\ge \ell$ and $|y-P|\le \fr{\ell_0}{2}$ satisfies $d(y)\le \fr{\ell}{2}$, then there exists a point $Q\in B(P,\ell_0)$ such that
\[
|x-y|^2\ge \fr{1}{16}(|y-Q|^2 + d(x)^2)\ge \fr{1}{16}(|y'-Q'|^2 + d(x)^2)
\]
Indeed, if $|x-P|\ge \ell_0$ we take $Q=P$ because then $|x-y| =|x-P+P-y|\ge \ell_0-\fr{\ell_0}{2}$, so $|x-y|\ge \fr{|y-Q|}{2}$. But also $|x-y|\ge\fr{d(x)}{2}$ because there exists a point $P_1 =(p, F(p))\in\pa\Omega$ such that $|y-P_1| = d(y)\le\fr{\ell}{2}$ while obviously $|x-P_1|\ge d(x)\ge\ell$.
If, on the other hand $|x-P|< \ell_0$, then $x$ is in the neighborhood of $P$ and we take $Q=x$. Because $y-P_1= (y'-p, y_d-F(p))$ we have
\[
d(y)\le |y_d-F(y')| \le \fr{11}{10}d(y)
\]
for $y\in B(P,\ell_0)$. We take a partition of unity of the form $1= \psi_0 +\sum_{j=1}^N\psi_j$ 
with $\psi_k\in C_0^{\infty}(\Rr^d)$, subordinated to the cover of the boundary with neighborhoods as above, and with $\psi_0$ supported in
$d(x)\ge \fr{\ell_0}{4}$, identically 1 for $d(x)\ge \fr{\ell_0}{2}$, $\psi_j$ supported near the boundary $\pa\Omega$ in balls of size $2\ell_0$ and identically 1 on balls of radius $\ell_0$.  

The cutoff will be taken of the form 
$\chi= \alpha_0 +\sum_{j=1}^N \chi_j(\fr{y_d-F(y')}{\ell})\alpha_j(y)$, where of course the meaning of $y$ changes in each neighborhood. The smooth functions $\chi_j(z)$, are identically zero for $|z|\le \fr{11}{40}$ and identically 1 for $|z| \ge\fr{10}{22}$. The integrals in (\ref{chij}) and (\ref{nachij}) reduce to integrals of the type 
\[
\ba
\int_{y_d>F(y'), |y'|\le\ell_0}\fr{\left(1-\chi_1\left(\fr{y_d-F(y')}{\ell}\right)\right)}{|x-y|^{d+j}}dy\le C\left(\int_{0}^{\infty}\left(1-\chi_1\left(\fr{u}{\ell}\right)\right)du\right)\left(\int_{\Rr^{d-1}}\fr{dy'}{\left({|y'-Q'|^2 +d(x)^2}\right)^{\fr{d+j}{2}}}\right)\\
\le
C\ell d(x)^{-1-j}\le Cd(x)^{-j}
\ea
\]
and
\[
\ba
\int_{y_d>F(y'), |y'|\le\ell_0}\fr{\left|\na_y\chi_1\left(\fr{y_d-F(y')}{\ell}\right)\right|}{|x-y|^{d-\alpha}}dy\le C\left(\int_{-{\infty}}^{\infty}|\na\chi_1(z)|dz\right)\left(\int_{\Rr^{d-1}}\fr{dy'}{\left({|y'-Q'|^2 +d(x)^2}\right)^{\fr{d-\alpha}{2}}}\right)\\
\le
C d(x)^{-(1-\alpha)}.
\ea
\]
This completes the proof.

We recall from (\cite{ci}) that the C\'{o}rdoba-C\'{o}rdoba inequality (\cite{cc}) holds in bounded domains. In fact, more is true: there is a lower bound that provides a strong boundary repulsive term:

\beg{prop}{\la{cordoba}} Let $\Omega$ be a bounded domain with smooth boundary. Let $0\le s<2$. There exists a constant $c>0$ depending only on the domain $\Omega$ and on $s$, such that, for any
$\Phi$, a $C^2$ convex function satisfying $\Phi(0)= 0$, and any  $f\in C_0^{\infty}(\Omega)$, the inequality
\be
\Phi'(f)\l^s f - \l^s(\Phi(f))\ge \fr{c}{d(x)^s}\left(f\Phi'(f)-\Phi(f)\right)
\la{cor}
 \ee
holds pointwise in $\Omega$.
\end{prop}
The proof follows in a straightforward manner from the proof of (\cite{ci}) using convexity, approximation, and the lower bound (\ref{lambdatheta}).
We prove below two nonlinear lower bounds for the case $\Phi(f)= \fr{f^2}{2}$, one when $f$ is a localized finite difference, and one when $f$ is a localized first derivative. The proof of Proposition \ref{cordoba} can be left as an exercise, following the same pattern as below.

\beg{thm}\la{nlmax}
Let $f\in L^{\infty}(\Omega)$ be smooth enough ($C^2$, e.g.) and vanish at the boundary, $f\in{\mathcal{D}}(\l^{s})$ with $0\le s<2$. 
Then 
\be
\ba
D(f) = f\l^{s} f -\fr{1}{2}\l^{s} f^2\\
 = \gamma_0\int_0^{\infty}t^{-1-\fr{s}{2}}dt\int_{\Omega}H_D(x,y,t)(f(x)-f(y))^2dy + \gamma_0 f^2(x)\int_0^{\infty}t^{-1-\fr{s}{2}}\left[1-e^{t\D}1\right](x)dt\\
= \gamma_0\int_0^{\infty}t^{-1-\fr{s}{2}}dt\int_{\Omega}H_D(x,y,t)(f(x)-f(y))^2dy + \; f^2(x)\fr{1}{2}\l^{s} 1.
\ea
\la{d}
\ee
holds for all $x\in\Omega$. Here $\gamma_0 =\fr{c_{s}}{2}$ with $c_s$ of (\ref{rep}).  Let $\ell>0$ be a small number and
let $\chi\in C_0^{\infty}(\Omega)$, $0\le\chi\le 1$ be a good cutoff function, with $\chi(y)=1$ for $d(y)\ge \fr{\ell}{2}$, $\chi(y) =0$ for $d(y)\le\fr{\ell}{4}$ and with 
$|\na\chi(y)|\le \fr{C}{\ell}$. 
There exist constants $\gamma_1>0$ and $M>0$ depending on $\Omega$ such that,
if $q(x)$ is a smooth function in $L^{\infty}(\Omega)$  then
if 
\[
f(x)=\chi(x)(\delta_h q(x)) =\chi(x)(q(x+h)-q(x))
\]
then
\be
D(f) = (f\l^{s} f)(x) -\fr{1}{2}(\l^{s} f^2)(x) \ge \gamma_1 |h|^{-s}\fr{|f_d(x)|^{2+s}}{\|q\|_{L^{\infty}}^{s}} + \gamma_1\fr{f^2(x)}{d(x)^{s}}
\la{nlb}
\ee
holds pointwise in $\Omega$ when $|h|\le\fr{\ell}{16}$, and $d(x)\ge \ell$ with
\be
|f_d(x)| = \left\{
\ba
|f(x)|,\quad\quad {\mbox{if}}\;\; |f(x)| \ge M\|q\|_{L^{\infty}(\Omega)}\fr{|h|}{d(x)},\\
0,\quad\quad\quad\quad {\mbox{if}}\;\; |f(x)| \le M\|q\|_{L^{\infty}(\Omega)}
\fr{|h|}{d(x)}.
\ea
\right.
\la{fd}
\ee
\end{thm}
\noindent{\bf{Proof.}} We start by proving (\ref{d}):
\[
\ba
f(x)\l^{s}f(x) - \fr{1}{2}\l^{s} f^2(x) \\
= c_{s}\int_0^{\infty}t^{-1-\fr{s}{2}}\int_{\Omega}\left\{\left[\fr{1}{|\Omega|} f(x)^2 - f(x)H_D(t,x,y)f(y)\right]- \fr{1}{2|\Omega|}f^{2}(x) + \fr{1}{2}H_D(t,x,y)f^2(y)\right\}dy\\
=c_{s}\int_0^{\infty}t^{-1-\fr{s}{2}}dt\int_{\Omega}\left\{\fr{1}{2}\left[H_D(t,x,y)(f(x)-f(y))^2\right]  + \fr{1}{2}f^2(x)\left[\fr{1}{|\Omega|} -H_D(t,x,y)\right]\right\}dy \\
= c_{s}\int_0^{\infty}t^{-1-\fr{s}{2}}dt\int_{\Omega}\left\{\fr{1}{2}\left[H_D(t,x,y)(f(x)-f(y))^2\right]dy + \fr{1}{2}f^2(x)\left[1-e^{t\D}1\right](x)\right\} \\
= c_{s}\int_0^{\infty}t^{-1-\fr{s}{2}}dt\int_{\Omega}\fr{1}{2}\left[H_D(t,x,y)(f(x)-f(y))^2\right]dy + \fr{1}{2}f^2(x)\l^{s}1.
\ea
\]
It follows that
\be
\ba
\left(f\l^{s}f - \fr{1}{2}\l^{s} f^2\right)(x) \\\ge \fr{1}{2}c_{s}\int_0^{\infty}\psi\left(\fr{t}{\tau}\right)t^{-1-\fr{s}{2}}dt\int_{\Omega}H_D(t,x,y)(f(x)-f(y))^2dy + \fr{1}{2}f^2(x)\l^{s}1
\ea
\la{lowone}
\ee
where $\tau>0$ is arbitrary and $0\le \psi(s)\le 1$ is a smooth function, vanishing identically for $0\le s\le 1$ and equal identically to $1$ for $s\ge 2$.
We restrict to $t\le T$,
\be
\ba
\left(f\l^{s}f - \fr{1}{2}\l^{s} f^2\right)(x)\\ \ge
\fr{1}{2}c_{s}\int_0^{T}\psi\left(\fr{t}{\tau}\right)t^{-1-\fr{s}{2}}dt\int_{\Omega}H_D(t,x,y)\left(f(x)-f(y)\right)^2dy +\fr{1}{2}f^2(x)\l^{s}1
\ea
\la{lowtwo}
\ee
and open brackets in (\ref{lowtwo}):
\be
\ba
\left(f\l^{s}f - \fr{1}{2}\l^{s} f^2\right)(x) 
\ge \fr{1}{2}f^2(x)c_{s}\int_0^T\psi\left(\fr{t}{\tau}\right)t^{-1-\fr{s}{2}}dt\int_{\Omega}H_D(t,x,y)dy\\ 
- f(x)c_{s}\int_0^T\psi\left(\fr{t}{\tau}\right)t^{-1-\fr{s}{2}}dt\int_{\Omega}H_D(t,x,y)f(y)dy + \fr{1}{2}f^2(x)\l^{s}1\\
\ge |f(x)|\left [ \fr{1}{2}|f(x)| I(x) - J(x)\right] +\fr{1}{2}f^2(x)\l^{s}1
\ea
\la{lowthree}
\ee
with
\be
I(x) = c_{s}\int_0^T\psi\left(\fr{t}{\tau}\right)t^{-1-\fr{s}{2}}dt\int_{\Omega}H_D(t,x,y)dy,
\la{ix}
\ee
and
\be
\ba
J(x) = c_{s}\left |\int_0^T\psi\left(\fr{t}{\tau}\right)t^{-1-\fr{s}{2}}dt\int_{\Omega}H_D(t,x,y)f(y)dy\right |\\
= c_{s}\left |\int_0^T\psi\left(\fr{t}{\tau}\right)t^{-1-\fr{s}{2}}dt\int_{\Omega}H_D(t,x,y)\chi(y)\delta_hq(y)dy\right |.
\ea
\la{jx}
\ee
We proceed with a lower bound on $I$ and an upper bound on $J$.  For the lower bound on $I$ we note that, in view of (\ref{thetalow}) and the fact that
\[
I(x) = c_{s}\int_0^T\psi\left(\fr{t}{\tau}\right)t^{-1-\fr{s}{2}}\Theta(x,t)dt
\]
we have
\[
\ba
I(x)\ge c_1\int_0^{\min(T, d^2(x))}\psi\left(\fr{t}{\tau}\right)t^{-1-\fr{s}{2}}dt\\
= c_1\tau^{-\fr{s}{2}}\int_1^{\tau^{-1}(\min(T, d^2(x)))}\psi(u)u^{-1-\fr{s}{2}}du.
\ea
\]
Therefore we have that
\be
I(x)\ge c_2 \tau^{-\fr{s}{2}}
\la{ilow}
\ee
with $c_2 = c_1\int_1^2\psi(u)u^{-1-\fr{s}{2}}du$, a positive constant depending only on $\Omega$ and $s$, provided $\tau$ is small enough,
\be
\tau \le \fr{1}{2}\min(T,  d^2(x)).
\la{taucond}
\ee
In order to bound $J$ from above we use (\ref{grup}) with $\alpha=0$. Now 
\[
\ba
J \le c_{s}\left |\int_0^T\psi\left(\fr{t}{\tau}\right)t^{-1-\fr{s}{2}}dt\int_{\Omega}\delta_{-h}H_D(t,x,y)\chi(y)q(y)dy\right | + \\
 c_{s}\left |\int_0^T\psi\left(\fr{t}{\tau}\right)t^{-1-\fr{s}{2}}dt\int_{\Omega}H_D(t,x,y-h)(\delta_{-h}\chi(y))q(y)dy\right |
\ea
\]
We have that 

\[
J_2 = c_{s}\left |\int_0^T\psi\left(\fr{t}{\tau}\right)t^{-1-\fr{s}{2}}dt\int_{\Omega}H_D(t,x,y-h)(\delta_{-h}\chi(y))q(y)dy\right |\le C_{6}\fr{|h|}{d(x)}\|q\|_{L^{\infty}}\tau^{-\fr{s}{2}}.
\]
Indeed, 
\[
t^{-{\fr{d}{2}}}e^{-\fr{|x-y|^2}{Kt}}\le C_K |x-y|^{-d}
\]
so the bound follows from (\ref{hb}) and (\ref{nachij}). On the other hand, 
\[
\ba
J_1=  c_{s}\left |\int_0^T\psi\left(\fr{t}{\tau}\right)t^{-1-\fr{s}{2}}dt\int_{\Omega}\delta_{-h}H_D(t,x,y)\chi (y)q(y)dy\right | \\
\le
\|q\|_{L^{\infty}(\Omega)}|h|\int_0^T\psi\left(\fr{t}{\tau}\right)t^{-1-\fr{s}{2}}dt\int_{\Omega}|\na_y H_D(t,x,y)|dy
\ea
\]
and therefore, in view of (\ref{grup})
\[
J_1\le C_1 |h|\|q\|_{L^{\infty}(\Omega)}\int_0^T\psi\left(\fr{t}{\tau}\right )t^{-\fr{3}{2}-\fr{s}{2}}dt
\]
and therefore
\be
J_1 \le C_7|h|\|q\|_{L^{\infty}(\Omega)}\tau^{-\fr{1}{2}-\fr{s}{2}}
\la{jupone}
\ee
with
\[
C_7 = C_1\int_1^{\infty}\psi(u)u^{-\fr{3}{2}-\fr{s}{2}}du
\]
a constant depending only on $\Omega$ and $s$. In conclusion
\be
|J| \le C_8\tau^{-\fr{s}{2}}|h| (\tau^{-\fr{1}{2}} + d(x)^{-1})\|q\|_{L^{\infty}}.
\la{jup}
\ee
Now, because of the lower bound (\ref{lowthree}), 
if we can choose $\tau$ so that
\[
J(x) \le \fr{1}{4} |f(x)|I(x)
\]
then it follows that
\be
\left[f\l^{s}f - \fr{1}{2}\l^{s} f^2\right](x) \ge \fr{1}{4}f^2(x)I(x) + \fr{1}{2}f^2(x)\l^{s}1.
\la{lowfour}
\ee
Because of the bounds (\ref{ilow}), (\ref{jup}), if 
\[
|f(x)|\ge \fr{8C_8}{c_2}\fr{|h|}{d(x)}\|q\|_{L^{\infty}},
\]
then a choice
\be
\tau(x)^{-\fr{1}{2}} = C_9{\|q\|_{L^{\infty}}^{-1}}|f(x)||h|^{-1}
\la{tauchoice}
\ee
with $C_9 = c_2 (8C_8)^{-1}$ achieves the desired bound. This concludes the proof.

We are providing now a lower bound for $D(f)$ for a different situation.

\beg{thm}\la{nlmaxd}
Let $\ell>0$ be a small number and let $\chi\in C_0^{\infty}(\Omega)$, $0\le\chi\le 1$ be a good cutoff function, with $\chi(y)=1$ for $d(y)\ge \fr{\ell}{2}$, $\chi(y) =0$ for $d(y)\le\fr{\ell}{4}$ and with $|\na\chi(y)|\le \fr{C}{\ell}$.
There exist constants $\gamma_2>0$ and $M>0$ depending on $\Omega$  such that, if $q(x)$ is a smooth function in $C^{\alpha}(\Omega)$ with $0<\alpha<1$ and 
\[
f(x)=\chi(x)\na q(x),
\]
then
\be
D(f) = (f\l^{s} f)(x) -\fr{1}{2}(\l^{s} f^2)(x) \ge \gamma_2 \fr{|f_d(x)|^{2+\fr{s}{1-\alpha}}}{\|q\|_{C^{\alpha}(\Omega)}^{\fr{s}{1-\alpha}}}(d(x))^{\fr{s\alpha}{1-\alpha}} + \gamma_1\fr{f^2(x)}{d(x)^{s}}
\la{nlbd}
\ee
holds pointwise in $\Omega$ when  $d(x)\ge \ell$, with
\be
|f_d(x)| = \left\{
\ba
|f(x)|,\quad\quad {\mbox{if}}\;\; |f(x)| \ge M\|q\|_{L^{\infty}(\Omega)}(d(x))^{-1},\\
0,\quad\quad\quad\quad {\mbox{if}}\;\; |f(x)| \le M\|q\|_{L^{\infty}(\Omega)}(d(x))^{-1}.
\ea
\right.
\la{fdd}
\ee
\end{thm}
\noindent{\bf Proof.} We follow exactly the proof of Theorem \ref{nlmax} up to, and including the definition of $I(x)$ given in (\ref{ix}). In particular, the lower bound (\ref{ilow}) is still valid, provided $\tau$ is small enough (\ref{taucond}). The term $J$ starts out the same, but is treated slightly differently,
\be
\ba 
J(x) = c_{s}\left |\int_0^T\psi\left(\fr{t}{\tau}\right)t^{-1-\fr{s}{2}}dt\int_{\Omega}H_D(t,x,y)f(y)dy\right |\\
= c_{s}\left |\int_0^T\psi\left(\fr{t}{\tau}\right)t^{-1-\fr{s}{2}}dt\int_{\Omega}H_D(t,x,y)\chi(y)\na_y (q(y)-q(x))dy\right |.
\ea
\la{jxd}
\ee
In order to bound $J$ we use (\ref{nachij}) and (\ref{grup}). 

\[
\ba
|J(x)| \le c_{s}\left |\int_0^T\psi\left(\fr{t}{\tau}\right)t^{-1-\fr{s}{2}}dt\int_{\Omega}\pa_yH_D(t,x,y)\chi(y)(q(y)-q(x))dy\right | + \\
 c_{s}\left |\int_0^T\psi\left(\fr{t}{\tau}\right)t^{-1-\fr{s}{2}}dt\int_{\Omega}H_D(t,x,y)(\na\chi(y))(q(y)-q(x))dy\right |\\
= J_1(x) + J_2(x) 
\ea
\]
We have from (\ref{hb}) and (\ref{nachij}), as before,
\[
J_2(x)  = c_{s}\left |\int_0^T\psi\left(\fr{t}{\tau}\right)t^{-1-\fr{s}{2}}dt\int_{\Omega}H_D(t,x,y)(\na\chi(y))(q(y)-q(x))dy\right |\le Cd(x)^{-1}\|q\|_{L^{\infty}}\tau^{-\fr{s}{2}}.
\]
 On the other hand,
\[
\ba
J_1(x) =  c_{s}\left |\int_0^T\psi\left(\fr{t}{\tau}\right)t^{-1-\fr{s}{2}}dt\int_{\Omega}\na_yH_D(t,x,y)\chi (y)(q(y)-q(x))dy\right | \\
\le
c_s(d(x))^{-\alpha}\|q\|_{C^{\alpha}(\Omega)}\int_0^T\psi\left(\fr{t}{\tau}\right)t^{-1-\fr{s}{2}}dt\int_{\Omega\cap |x-y|\le d(x)}|\na_y H_D(t,x,y)||x-y|^{\alpha}dy \\
+ c_s\|q\|_{L^{\infty}}\int_0^T\psi\left(\fr{t}{\tau}\right)t^{-1-\fr{s}{2}}dt\int_{\Omega\cap |x-y|\ge d(x) }|\na_y H_D(t,x,y)|dy \\
= J_{11}(x) + J_{12}(x).
\ea
\]
In view of (\ref{grup})
\[
J_{11}(x)\le C_1 d(x)^{-\alpha} \|q\|_{C^{\alpha}(\Omega)}\int_0^T\psi\left(\fr{t}{\tau}\right )t^{-\fr{3-\alpha}{2}-\fr{s}{2}}dt
\]
and so
\be
J_{11}(x) \le C_2(d(x))^{-\alpha}\|q\|_{C^{\alpha}(\Omega)}\tau^{-\fr{1-\alpha}{2}-\fr{s}{2}}
\la{juponed}
\ee
with
\[
C_2 = C_1\int_1^{\infty}\psi(z)z^{-\fr{3-\alpha}{2}-\fr{s}{2}}dz
\]
a constant depending only on $\Omega$ and $s$. Regarding $J_{12}(x)$ we have in view of (\ref{grby}) 
\[
J_{12}(x) \le C\|q\|_{L^{\infty}(\Omega)}\int_0^T\psi\left(\fr{t}{\tau}\right )t^{-1-\fr{s}{2}}\left (\fr{1}{\sqrt{t}} + \fr{1}{d(x)}\right )e^{-\fr{d(x)^2}{2Kt}}dt\le C_3\tau^{-\fr{s}{2}}d(x)^{-1}\|q\|_{L^{\infty}(\Omega)}
\]
because, in view of (\ref{phione}) 
\[
\fr{w_1(y)}{|x-y|}\le C_0\fr{d(y)}{|x-y|} \le C_0\fr{d(y)}{d(x)}
\]
on the domain of integration. 

In conclusion
\be
|J(x)| \le C_3\tau^{-\fr{s}{2}}\left[\tau^{-\fr{1-\alpha}{2}}(d(x))^{-{\alpha}}\|q\|_{C^{\alpha}(\Omega)} + d(x)^{-1}\|q\|_{L^{\infty}(\Omega)}\right].
\la{jupd}
\ee
The rest is the same as in the proof of Theorem \ref{nlmax}:
If $|f(x)|\ge Md(x)^{-1}\|q\|_{L^{\infty}(\Omega)}$ for suitable $M$, ($M= 8C_3c_2^{-1}$) then we choose $\tau$ such that
\[
\fr{|f(x)|}{\|q\|_{C^{\alpha}(\Omega)}} = M\tau^{-\fr{1-\alpha}{2}}(d(x))^{-\alpha},
\]
and this yields  $|f(x)| I\ge 4 |J(x)|$, and consequently, in view of
(\ref{lowfour}) which is then valid, the result (\ref{nlbd}) is proved.

We specialize from now on to $s=1$. 

\section{Bounds for Riesz transforms}
We consider $u$ given in (\ref{u}),  
\[
u = \na^{\perp}\l^{-1}\theta.
\]
We are interested in estimates of $u$ in terms of $\theta$, and in particular estimates of finite differences and the gradient.
We fix a length scale $\ell$ and take a good  cutoff function $\chi\in C_0^{\infty}(\Omega)$ which satisfies $\chi(x) =1 $ if $d(x)\ge \fr{\ell}{2}$, $\chi(x) = 0$ if $d(x)\le \fr{\ell}{4}$, $|\na\chi(x)|\le C\ell^{-1}$, (\ref{chij}) and (\ref{nachij}). 
We take $|h|\le \fr{\ell}{14}$. In view of the representation
\be
\l^{-1} = c\int_0^{\infty}t^{-\fr{1}{2}}e^{t\D}dt
\la{lambdaminusone}
\ee
we have on the support of $\chi$
\be
\delta_h u(x) = c\int_0^{\infty}t^{-\fr{1}{2}}dt\int_{\Omega}\delta_h^x\na_x^{\perp}H_D(x,y,t)\theta(y)dy.
\la{duh}
\ee
We split
\be
\delta_h u = \delta_h u^{in} + \delta_h u^{out}
\la{split}
\ee
with
\be
\delta_h u(x)^{in} = c\int_0^{\rho^2}t^{-\fr{1}{2}}dt\int_{\Omega}\delta_h^x\na_x^{\perp}H_D(x,y,t)\theta(y)dy
\la{duhin}
\ee
and $\rho=\rho(x,h)>0$ a length scale to be chosen later but it will be smaller than the distance from $x$ to the boundary of $\Omega$:
\be
\rho \le c d(x).
\la{rholess}
\ee
We represent 
\be
\delta_h u^{in}(x) = u_h(x) + v_h(x)
\la{deltahu}
\ee
where
\be 
u_h(x) = c\int_0^{\rho^2}t^{-\fr{1}{2}}dt\int_{\Omega}\na_x^{\perp}H(x,y,t)(\chi(y)\delta_h\theta(y)-\chi(x)\delta_h\theta(x))dy 
\la{uh}
\ee
and where
\be
v_h(x) = e_1(x) + e_2(x) + e_3(x) + \chi(x)\delta_h\theta(x)e_4(x)
\la{vh}
\ee
with
\be
e_1(x) = c\int_0^{\rho^2}t^{-\fr{1}{2}}dt\int_{\Omega}\na_x^{\perp}(H_D(x+h,y,t)-H_D(x,y,t))(1-\chi(y))\theta(y)dy,
\la{e1}
\ee
\be
e_2(x) =  c\int_0^{\rho^2}t^{-\fr{1}{2}}dt\int_{\Omega}\na_x^{\perp}(H_D(x+h,y,t)-H_D(x,y-h,t))\chi(y)\theta(y)dy,
\la{e2}
\ee
\be
e_3(x) =c\int_0^{\rho^2}t^{-\fr{1}{2}}dt\int_{\Omega}\na_x^{\perp}H_D(x,y,t)(\chi(y+h)-\chi(y))\theta(y+h)dy,
\la{e3}
\ee
and
\be
e_4(x) =c\int_0^{\rho^2}t^{-\fr{1}{2}}dt\int_{\Omega}\na_x^{\perp}H_D(x,y,t)dy.
\la{e4}
\ee
We used here the fact that $(\chi\theta)(\cdot)$ and  $(\chi\theta)(\cdot + h)$
are compactly supported in $\Omega$ and hence 
\[
\int_{\Omega}\na_x^{\perp}H_D(x,y-h,t)\chi(y)\theta(y)dy = \int_{\Omega}\na_x^{\perp}H_D(x,y,t)\chi(y+h)\theta(y+h)dy.
\]
The following elementary lemma will be used in several instances:
\beg{lemma}\la{intpk}
Let $\rho>0$, $p>0$. Then
\be
\int_0^{\rho^2} t^{-1-\fr{m}{2}}\left(\fr{p}{\sqrt{t}}\right)^je^{-\fr{p^2}{Kt}}dt \le  C_{K,m,j}p^{-m}
\la{pbeta}
\ee
if $m\ge 0$, $j\ge 0$, $m+j>0$, and
\be
\int_0^{\rho^2} t^{-1}e^{-\fr{p^2}{Kt}}dt \le  C_K\left(1+ 2\log_{+}\left(\fr{\sqrt{K}\rho}{p}\right)\right)
\la{pzero}
\ee
if $m=0$ and $j=0$, with constants $C_{K,m,j}$ and $C_K$ independent of $\rho$ and $p$. Note that when $m+j>0$, $\rho=\infty$ is allowed.

\end{lemma}

We start estimating the terms in (\ref{vh}).
For $e_1$ we use the inequality (\ref{naxnaxb}), and it then follows from Lemma \ref{intpk} with $m= {d} +1$ that
\[
|e_1(x)|\le C|h|\|\theta\|_{L^{\infty}}\int_0^1d\lambda \int_{\Omega}\fr{1}{|x +\lambda h -y |^{d+1}}(1-\chi(y))dy
\]
and therefore we have from (\ref{chij}) that
\be
|e_1(x)| \le C \|\theta\|_{L^{\infty}} \fr{|h|}{d(x)}
\la{e1b}
\ee
holds for $d(x)\ge \ell$. Concerning $e_3$ we use Lemma (\ref{intpk}) with
$m= d$ and $j=0,1$ in conjunction with (\ref{grbx}) and obtain
\[
|e_3(x)|\le C{|h|}\|\theta\|_{L^{\infty}}\int_{\Omega}|\na\chi(y)|\fr{1}{|x-y|^d}dy
\]
and therefore we obtain from (\ref{nachij})
\be
|e_3(x)| \le C \|\theta\|_{L^{\infty}}\fr{|h|}{d(x)}
\la{e3b}
\ee
holds for $d(x)\ge \ell$. Regarding $e_4$ we can split it into
\[
e_4 (x) = e_5(x) + e_6(x)
\]
with 
\[
e_{5}(x) = \int_0^{\rho^2}t^{-\fr{1}{2}}\int_{\Omega}\na_x^{\perp}H_D(x,y,t)\chi(y)dy
\]
and 
\[
e_{6}(x) = \int_0^{\rho^2}t^{-\fr{1}{2}}\int_{\Omega}\na_x^{\perp}H_D(x,y,t)(1-\chi(y))dy.
\]
Now $e_6$ is bounded using the Lemma (\ref{intpk}) with $m= d$ and $j=0,1$ in conjunction with (\ref{grbx}) and (\ref{chij}) and obtain
\be
|e_6(x)|\le C\int_{\Omega}\fr{(1-\chi(y))}{|x-y|^d}dy \le C
\la{e6b}
\ee
for $d(x)\ge \ell$, with a constant independent of $\ell$. For $e_5$ we use the fact that $\chi $ is a fixed smooth function which vanishes at the boundary.

In order to bound the terms $e_2$ and $e_5$ we need to use additional information, namely the inequalities (\ref{cancel1}) and (\ref{cancel2}). For $e_5$ we write
\[
\ba
e_5(x) = \int_0^{\rho^2}t^{-\fr{1}{2}}dt\int_{\Omega}\left(\na_x^{\perp}H_D(x,y,t) + \na_y^{\perp} H_D(x,y,t)\right)\chi(y)dy \\
+ \int_0^{\rho^2}t^{-\fr{1}{2}}dt\int_{\Omega}H_D(x,y,t)\na_y^{\perp}\chi(y)dy,
\ea
\]
and using (\ref{cancel1}) and Lemma \ref{intpk} with $m=0$, $j=0$ and (\ref{nachij}) we obtain the bound
\[
|e_5(x)|\le C\left(1+ \log_{+}\left(\fr{\rho}{d(x)}\right)\right) + C\rho\int_{\Omega}\fr{|\na \chi(y)|}{|x-y|^d}dy
\]
and  therefore, in view of (\ref{nachij}) and $\rho\le d(x)$ we have
\be
|e_5(x)|\le C
\la{e5b}
\ee
for $d(x)\ge \ell$, with $C$ depending on $\Omega$ but not on $\ell$.
Consequently, we have
\be
|e_4(x)|\le C
\la{e4b}
\ee
for  $d(x)\le \ell$, with a constant $C$ depending on $\Omega$ only.
In order to estimate $e_2$ we write
\be
H_D(x+h,y,t)- H_D(x,y-h,t) = h\cdot\int_0^1 (\na_x+\na_y)H_D(x+\lambda h, y + (\lambda-1)h,t)d\lambda
\la{transh}
\ee
and use (\ref{cancel2}) and Lemma \ref{intpk} with $m=1$, $j=0$ to obtain
\[
\ba
|e_2(x)|\le |h|\|\theta\|_{L^{\infty}}\int_0^1d\lambda \int_0^{\rho^2}t^{-\fr{1}{2}}dt\int_{\Omega}|\na_x^{\perp}(\na_x +\na_y)H_D(x+\lambda h, y+(\lambda-1)h)||\chi(y)|dy\\
\le C|h|\|\theta\|_{L^{\infty}}\int_0^1d\lambda\int_0^{\rho^2}t^{-\fr{3}{2}}e^{-\fr{d(x)^2}{4Kt}}dt
\ea
\]
and thus
\be
|e_2(x)|\le C\|\theta\|_{L^{\infty}}\fr{|h|}{d(x)}
\la{e2b}
\ee
holds for $d(x)\ge \ell$.
Summarizing, we have that
\be
|v_h(x)| \le C\|\theta\|_{L^{\infty}}\fr{|h|}{d(x)} + C|\delta_h\theta(x)|
\la{vhbo}
\ee
for $d(x)\ge \ell$. We now estimate $u_h$ using (\ref{grbx}) and a Schwartz inequality
\[
\ba 
|u_h(x)| \le c\int_0^{\rho^2}t^{-1}\int_{\Omega}\left(1+\fr{|x-y|}{\sqrt{t}}\right)H_D(x,y,t)(\chi(\delta_h\theta)(y)-\chi\delta_h\theta(x))dy\\
\le \sqrt{\rho}\left\{\int_0^{\rho^2}t^{-\fr{3}{2}}dt\int_{\Omega}H_D(x,y,t)(\chi(\delta_h\theta)(y)-\chi\delta_h\theta)^2dy\right\}^{\fr{1}{2}}
\ea
\]
We have therefore
\be
|u_h(x)|\le C\sqrt{\rho D(f)(x)}.
\la{uhbo}
\ee
where $f=\chi\delta_h\theta$ and $D(f)$ is given in Theorem \ref{nlmax}.
Regarding $\delta_h u^{out}$ we have 
\be
|\delta_h u^{out}(x)| \le C\|\theta\|_{L^{\infty}}\fr{|h|}{\rho}
\la{dhoutb}
\ee
in view of (\ref{naxnaxb}). Putting together the estimates (\ref{vhbo}), (\ref{uhbo}) and (\ref{dhoutb}) we have

\beg{prop} Let $\chi$ be a good cutoff, and let $u$ be defined by 
(\ref{u}). Then
\be
|\delta_h u(x)| \le  C\left(\sqrt{\rho D(f)(x)} + \|\theta\|_{L^{\infty}}\left(\fr{|h|}{d(x)}+ \fr{|h|}{\rho}\right) + |\delta_h\theta(x)|\right)
\la{dhub}
\ee
holds for $d(x)\ge\ell$, $\rho\le cd(x)$, $f=\chi\delta_h \theta$ and  with $C$ a constant depending on $\Omega$.
\end{prop}

Now we will obtain similar estimates for $\na u$. We start with the representation
\be
\na u(x) = \na u^{in}(x) + \na u^{out}(x)
\la{nauinout}
\ee
where
\be
\na u^{in}(x) = c\int_0^{\rho^2}t^{-\fr{1}{2}}\int_{\Omega}\na_x\na_x^{\perp}H_D(x,y,t)\theta(y)dy
\la{nauin}
\ee
and $\rho= \rho(x) \le c d(x)$.
In view of (\ref{naxnaxb}) we have
\be
|\na u^{out}(x)|\le \fr{C}{\rho}\|\theta\|_{L^{\infty}(\Omega)}
\la{nauoutb}
\ee
We split now 
\be
\na u^{in}(x) = g(x) +  g_1(x) + g_2(x) + g_3(x) + g_4(x)f(x)
\la{nasplit}
\ee
where 
\be
f(x) = \chi(x)\na\theta(x)
\la{fnatheta}
\ee
and with
\be
g(x) =  c\int_0^{\rho^2}t^{-\fr{1}{2}}\int_{\Omega}\na_x^{\perp}H_D(x,y,t)(f(y)-f(x))dy,
\la{g}
\ee
and 
\be
g_1(x) = c\int_0^{\rho^2}t^{-\fr{1}{2}}\int_{\Omega}\na_x\na_x^{\perp}(H_D(x,y,t)(1-\chi(y))\theta(y)dy,
\la{g1}
\ee
\be
g_2(x) =  c\int_0^{\rho^2}t^{-\fr{1}{2}}\int_{\Omega}\na_x^{\perp}(\na_x+ \na_y)H_D(x,y,t)\chi(y)\theta(y)dy,
\la{g2}
\ee
\be
g_3(x) =c\int_0^{\rho^2}t^{-\fr{1}{2}}\int_{\Omega}\na_x^{\perp}H_D(x,y,t)(\na_y\chi(y))\theta(y)dy,
\la{g3}
\ee
and
\be
g_4(x) =c\int_0^{\rho^2}t^{-\fr{1}{2}}\int_{\Omega}\na_x^{\perp}H_D(x,y,t)dy.
\la{g4}
\ee
Now
\be
|g_1(x)| \le \fr{C}{d(x)}\|\theta\|_{L^{\infty}(\Omega)}
\la{g1b}
\ee
holds for $d(x)\ge \ell$ because of (\ref{naxnaxb}), time integration using  Lemma \ref{intpk} and then use of (\ref{chij}).
For $g_2(x)$ we use (\ref{cancel2}) and then Lemma \ref{intpk} to obtain
\be
|g_2(x)| \le \fr{C}{d(x)}\|\theta\|_{L^{\infty}(\Omega)}
\la{g2b}
\ee
for $d(x)\ge \ell$. Now
\be
|g_3(x)| \le \fr{C}{d(x)}\|\theta\|_{L^{\infty}(\Omega)}
\la{g3b}
\ee
holds because of (\ref{grbx}), Lemma \ref{intpk} and then use of (\ref{nachij}).
Regarding $g_4$, in view of
\be
\int_{\Omega}\na_y^{\perp}H_D(x,y,t)dy = 0
\la{intnah}
\ee
we have
\[
g_4(x) = c\int_0^{\rho^2}t^{-\fr{1}{2}}\int_{\Omega}\left(\na_x^{\perp} + \na_y^{\perp}\right)H_D(x,y,t)dy
\]
and, we thus obtain from (\ref{cancel1}) and from Lemma \ref{intpk} with $m=j=0$
\be
|g_4(x)| \le C
\la{g4b}
\ee
because $\rho\le cd(x)$. 

Finally we have using a Schwartz inequality like for (\ref{uhbo})
\be
|g(x)| \le C\sqrt{\rho D(f)}.
\la{gb}
\ee
Gathering the bounds we have proved
\beg{prop}\la{naub}
Let $\chi$ be a good cutoff with scale $\ell$ and let $u$ be given by (\ref{u}). Then
\be
|\na u(x)| \le C\left(\sqrt{\rho D(f)} + \|\theta\|_{L^{\infty}(\Omega)}\left(\fr{1}{d(x)} + \fr{1}{\rho}\right) + |\na\theta(x)|\right)
\la{nauxb}
\ee
holds for $d(x)\ge \ell$, $\rho\le cd(x)$  and $f=\chi\na\theta$ with a constant $C$ depending on $\Omega$.
\end{prop}

\section{Commutators}
We consider the finite difference
\be
(\delta_h\l\theta)(x) = \l\theta(x+h)-\l\theta(x)
\la{dht}
\ee
with $d(x)\ge \ell$ and $|h|\le \fr{\ell}{16}$. We use a good cutoff $\chi$ again and denote
\be
f(x) = \chi(x)\delta_h\theta(x).
\la{f}
\ee
We start by computing
\be
\ba
(\delta_h\l\theta)(x) = (\l f)(x) + c\int_0^{\infty}t^{-\fr{3}{2}}dt\int_{\Omega}(H_D(x,y,t)-H_D(x+h,y,t))(1-\chi(y))\theta(y)dy\\
-c\int_0^{\infty}t^{-\fr{3}{2}}dt\int_{\Omega}(H_D(x+h,y,t)-H_D(x,y-h,t))\chi(y)\theta(y)dy\\
-c\int_0^{\infty}t^{-\fr{3}{2}}dt\int_{\Omega}H_D(x,y,t)(\delta_h\chi)(y)\theta(y+h)dy\\
 = (\l f)(x) + E_1(x) + E_2(x) + E_3(x). 
\ea
\la{compucom}
\ee

\beg{lemma}\la{commuh}
There exists a constant $\Gamma_0$ such that the commutator
\be
C_h(\theta) = \delta_h\l\theta -\l(\chi\delta_h\theta)
\la{comh}
\ee
obeys
\be
\left |C_h(\theta)(x)\right| \le \Gamma_0\fr{|h|}{d(x)^2}\|\theta\|_{L^{\infty}(\Omega)}
\la{commhb}
\ee
for $d(x)\ge \ell$, $|h|\le\fr{\ell}{16}$, $f=\chi\delta_h\theta$ and $\theta\in H_0^1(\Omega)\cap L^{\infty}(\Omega)$.
\end{lemma}
\noindent{\bf Proof.} We use (\ref{compucom}).  For $E_1(x)$ we use a similar argument as for $e_1$ leading to (\ref{e1b}), namely the inequality (\ref{hb}) and Lemma \ref{intpk} with $m=d+2$, $j=0$, and (\ref{chij}) to obtain
\[
|E_1(x)|\le C\fr{|h|}{d(x)^2}\|\theta\|_{L^{\infty}}.
\]
For $E_2$ we proceed in a manner analagous to the one leading to the bound (\ref{e2b}), by using (\ref{transh}), (\ref{cancel1}),  Lemma \ref{intpk} with $m=d+2$, $j=0$, and (\ref{nachij}) to obtain
\[
|E_2(x)|\le  C\fr{|h|}{d(x)^2}\|\theta\|_{L^{\infty}}.
\]

For $E_3$ we use 
\[
|E_3(x)|\le |h|\|\theta\|_{L^{\infty}}\int_0^{\infty}t^{-\fr{3}{2}}dt\int_{\Omega}H_D(x,y,t)|\na(\chi)(y)|dy
\]
and using Lemma \ref{intpk} with $m=d+1$, $j=0$ and (\ref{nachij}) we obtain
\[
|E_3(x)|\le  C\fr{|h|}{d(x)^2}\|\theta\|_{L^{\infty}},
\]
concluding the proof.

We consider now the commutator $[\na, \l]$.

\beg{lemma}\la{commd} There exists a constant $\Gamma_3$ depending on $\Omega$ such that for any smooth function $f$ vanishing at $\pa\Omega$ and any $x\in\Omega$ we have
\be
\left| [\na, \l]f(x)\right| \le \fr{\Gamma_3}{d(x)^2}\|f\|_{L^{\infty}(\Omega)}.
\la{nac}
\ee
If $\chi$ is a good cutoff with scale $\ell$ and if $\theta$ is a smooth bounded function in ${\mathcal{D}}\left(\l\right)$, then 
\be
C_{\chi}(\theta) = \na\l\theta -\l\chi\na\theta
\la{cchi}
\ee
obeys
\be
|C_{\chi}(\theta)(x)| = \left|(\na\l\theta - \l(\chi\na\theta))(x)\right| \le \fr{\Gamma_3}{d(x)^2}\|\theta\|_{L^{\infty}(\Omega)}
\la{nachi}
\ee
for $d(x)\ge \ell$, with a constant $\Gamma_3$ independent of $\ell$.
\end{lemma}
\noindent{\bf Proof.} We note that
\be
[\na,\l]f(x) = -c_1\int_0^{\infty}t^{-\fr{3}{2}}\int_{\Omega}\left(\na_xH_D(x,y,t)f(y) - H_D(x,y,t)\na_y f(y)\right)dy
\la{commde}
\ee
and therefore
\be
[\na,\l]f(x) = -c_1\int_0^{\infty}t^{-\fr{3}{2}}\int_{\Omega}\left(\na_x+\na_y\right)H_D(x,y,t)f(y)dy.
\la{commdef}
\ee
The inequality (\ref{nac}) follows from (\ref{cancel1}) and Lemma \ref{intpk}.
For the inequality (\ref{nachi}) we need also to estimate
\[
C(x) = c_s\left|\int_0^{\infty}t^{-\fr{3}{2}}\int_{\Omega}H_D(x,y,t)(\na\chi(y))\theta(y)dy\right|
\]
by the right hand side of (\ref{nachi}), and this follows from (\ref{nachij}) in view of (\ref{hb}).

\section{SQG: H\"{o}lder bounds}
We consider the equation (\ref{sqg}) with $u$ given by (\ref{u}) and with smooth initial data $\theta_0$ compactly supported  in $\Omega$.  We note that by the C\'{o}rdoba-C\'{o}rdoba inequality we have
\be
\|\theta(t)\|_{L^{\infty}}\le \|\theta_0\|_{L^{\infty}}.
\la{linf}
\ee
We prove the following uniform interior H\"{o}lder bound: 
\beg{thm}\la{glh} Let $\theta(x,t)$ be a smooth solution of (\ref{sqg}) in the smooth bounded domain $\Omega$. There exists a constant $0<\alpha<1$ depending 
only on $\|\theta_0\|_{L^{\infty}(\Omega)}$, and a constant $\Gamma>0$ depending on the domain $\Omega$ such that, for any $\ell>0$ sufficiently small
\be
\sup_{d(x)\ge \ell, \; |h|\le\fr{\ell}{16},\; t\ge 0}\fr{|\theta(x+h,t)-\theta(x,t)|}{|h|^{\alpha}} \le \|\theta_0\|_{C^{\alpha}} + \Gamma \ell^{-\alpha}\|\theta_0\|_{L^{\infty}(\Omega)}
\la{hldrbnd}
\ee
holds.
\end{thm}

\noindent{\bf{Proof.}} We take a good cutoff $\chi$ used above,  $|h|\le\fr{\ell}{16}$  and observe that, from the SQG equation we obtain the equation
\be
(\pa_t + u\cdot\na + (\delta_h u)\cdot\na_h)(\delta_h\theta)+ \l(\chi \delta_h \theta) + C_h(\theta) = 0
\la{le}
\ee
where $C_h(\theta)$ is the commutator given above in (\ref{comh}). Denoting (as before in (\ref{f})) $f=\chi\delta_h\theta$ we have after multiplying by $\delta_h\theta$ and using the fact that $\chi(x)=1$ for $d(x)\ge \ell$,
\be
\fr{1}{2}L_{\chi}\left(\delta_h\theta\right)^2 +  D(f) + (\delta_h\theta) C_h(\theta) = 0
\la{lchieq}
\ee
where
\be
L_{\chi}g = \pa_t g + u\cdot\na_x g + \delta_h u\cdot\na_h g  + \l(\chi^2 g) 
\la{lchi}
\ee
and $D(f)$ is given in Theorem \ref{nlmax}.

Multiplying by $|h|^{-2\alpha}$ where $\alpha>0$ will be chosen small to be small enough we obtain
\be
\fr{1}{2}L_{\chi}\left(\fr{\delta_h\theta(x)^2}{|h|^{2\alpha}}\right) + |h|^{-2\alpha} D(f) \le 2\alpha\fr{|\delta_h u|}{|h|}\left(\fr{\delta_h\theta(x)^2}{|h|^{2\alpha}}\right) + |C_h(\theta)| |\delta_h\theta||h|^{-2\alpha}.
\la{bu}
\ee
The factor $2\alpha$ comes from the differentiation $\delta_h u\cdot\na_h (|h|^{-2\alpha})$ and its smallness will be crucial below. 
Let us record here the inequality (\ref{d}) in the present case:
\be
D(f) \ge \gamma_1 |h|^{-1}\|\theta\|_{L^{\infty}}^{-1}|(\delta_h\theta)_d|^3 + \gamma_1 (d(x))^{-1} |\delta_h\theta|^2,
\la{d1}
\ee
valid pointwise, when $|h|\le \fr{\ell}{16}$ and $d(x)\ge\ell$, where
\[
|(\delta_h\theta)_d| = |\delta_h\theta|, \quad {\mbox{if}}\; |\delta_h\theta(x)|\ge M\|\theta\|_{L^{\infty}}\fr{|h|}{d(x)},
\]
and $|(\delta_h\theta)_d|=0$ otherwise.

We use now the estimates (\ref{dhub}), (\ref{commhb}) and a Young inequality for the term involving $\sqrt{\rho D(f)}$ to obtain
\be
\ba
\fr{1}{2}L_{\chi}\left(\fr{\delta_h\theta(x)^2}{|h|^{2\alpha}}\right) + \fr{1}{2}|h|^{-2\alpha} D(f) \le 
 C_1\alpha^2|h|^{-2-2\alpha}\rho |\delta_h\theta|^4 \\+ 
C_1\alpha \|\theta\|_{L^{\infty}}\left(\fr{1}{d(x)} + \fr{1}{\rho}\right)|h|^{-2\alpha}|\delta_h\theta|^2 + C_1\alpha|\delta_h\theta||h|^{-1-2\alpha}|\delta_h\theta|^2\\
+ \Gamma_0\fr{|h|}{d(x)^2}\|\theta\|_{L^{\infty}} |\delta_h\theta||h|^{-2\alpha}
\ea
\la{bo}
\ee
for $d(x)\ge \ell$, $|h|\le \fr{\ell}{16}$. Let us choose $\rho$ now. We set
\be
\rho = \left\{
\ba 
|\delta_h\theta(x)|^{-1}|h|\|\theta\|_{L^{\infty}}, \quad {\mbox{if}}\;
|\delta_h\theta(x)|\ge M_1\|\theta\|_{L^{\infty}}\fr{|h|}{d(x)},\\

d(x), \quad\quad {\mbox{if}}\quad |\delta_h\theta(x)|\le M_1\|\theta\|_{L^{\infty}}\fr{|h|}{d(x)},
\ea
\right.
\la{rhoc}
\ee
where we put
\be
M_1 = M + \sqrt{\fr{8\Gamma_0}{\gamma_1}} + 1,
\la{mone}
\ee
where $M$ is the constant from Theorem \ref{nlmax}, $\Gamma_0$ is the constant from (\ref{commhb}) and $\gamma_1$ is the constant from (\ref{d1}). This choice was made in order to use the lower bound on $D(f)$ to estimate the contribution due to the inner piece $u_h$ (see (\ref{uh})) of $\delta_h u$ and the contribution from the commutator $C_h(\theta)$. We distinguish two cases. The first case is when
$|\delta_h\theta(x)|\ge M_1\|\theta\|_{L^{\infty}}\fr{|h|}{d(x)}$. Then we have
\be
\ba
\fr{1}{2}L_{\chi}\left(\fr{\delta_h\theta(x)^2}{|h|^{2\alpha}}\right) + \fr{1}{2}|h|^{-2\alpha}D(f) \le  
 C_1\left[(\alpha\|\theta\|_{L^{\infty}})^2 +(2+\fr{1}{M_1})\alpha\|\theta\|_{L^{\infty}}\right]|\delta_h\theta|^3|h|^{-1-2\alpha}\|\theta\|_{L^{\infty}}^{-1} \\
+ \Gamma_0\fr{|h|}{d(x)^2}\|\theta\|_{L^{\infty}} |\delta_h\theta||h|^{-2\alpha}.
\ea
\la{boup} 
\ee
The choice of $M_1$ was such that, in this case
\[
\Gamma_0\fr{|h|}{d(x)^2}\|\theta\|_{L^{\infty}} |\delta_h\theta(x)||h|^{-2\alpha}
\le \fr{\gamma_1}{8} |\delta_h\theta(x)|^3|h|^{-1-2\alpha}\|\theta\|_{L^{\infty}}^{-1} .
\]
We choose now $\alpha$ by requiring
\be
\epsilon= \alpha\|\theta\|_{L^{\infty}}
\la{epsilong}
\ee
to satisfy
\be
C_1M_1^2(\epsilon^2 + (2+M_1^{-1})\epsilon) \le\fr{\gamma_1}{8}
\la{epsilongre}
\ee
and obtain from (\ref{boup})
\be
\fr{1}{2}L_{\chi}\left(\fr{|\delta_h\theta(x)|^2}{|h|^{2\alpha}}\right) + \fr{1}{4}|h|^{-2\alpha}D(f) \le 0
\la{bonu}
\ee
for $d(x)\ge\ell$, $|h|\le\fr{\ell}{16}$, in the case $|f|\ge M_1\|\theta\|_{L^{\infty}}\fr{|h|}{d(x)}$. 

The second case is when the opposite inequality holds, i.e, when
$|\delta_h\theta(x)|\le M_1\|\theta\|_{L^{\infty}}\fr{|h|}{d(x)}$. Then, using $\rho = d(x)$ we obtain from (\ref{bo})
\be
\ba
\fr{1}{2}L_{\chi}\left(\fr{\delta_h\theta(x)^2}{|h|^{2\alpha}}\right) + \fr{1}{2}|h|^{-2\alpha}D(f) \le  C_1(M_1^2\epsilon^2 + (M_1+2)\epsilon)\fr{1}{d(x)} (\delta_h\theta(x))^2|h|^{-2\alpha}\\ +
\Gamma_0d(x)^{-2}\|\theta\|_{L^{\infty}}|\delta_h\theta||h|^{1-2\alpha}\\
\le \fr{\gamma_1}{8d(x)}\left(\fr{\delta_h\theta(x)^2}{|h|^{2\alpha}}\right)
+ 2\Gamma_0M_1\|\theta\|_{L^{\infty}}^2d(x)^{-3}|h|^{2-2\alpha}.
\ea
\la{badu}
\ee
Summarizing, in view of the inequalities (\ref{bonu}) and (\ref{badu}), the damping term $\fr{\gamma_1}{d(x)}|\delta_h\theta(x)|^2$ in (\ref{d1}) and the choice of small $\epsilon$ in (\ref{epsilongre}), we have that
\be
L_{\chi}\left(\fr{\delta_h\theta(x)^2}{|h|^{2\alpha}}\right) + \fr{\gamma_1}{4d(x)}\left(\fr{\delta_h\theta(x)^2}{|h|^{2\alpha}}\right) \le
B
\la{fin}
\ee
holds for $d(x)\ge \ell$ and $|h|\le\fr{\ell}{16}$ where 

\be
B=  2(16)^{-2+2\alpha}\Gamma_0M_1\|\theta\|_{L^{\infty}}^2d(x)^{-1-2\alpha} =
\Gamma_1\fr{\gamma_1}{4}\|\theta\|_{L^{\infty}}^2d(x)^{-1-2\alpha}
\la{B}
\ee
with $\Gamma_1$ depending on $\Omega$.  Without loss of generality we may 
take $\Gamma_1> 4(16)^{2\alpha}$ so that
\[
\fr{|\delta_h\theta|^2}{|h|^{2\alpha}} <\Gamma_1\ell^{-2\alpha}\|\theta_0\|_{L^{\infty}}^2
\]
when $|h|\ge \fr{\ell}{16}$.  We note that
\be
L_\chi \left(\fr{\delta_h\theta(x)^2}{|h|^{2\alpha}}\right) +
 \fr{\gamma_1}{4d(x)}\left(\fr{\delta_h\theta(x)^2}{|h|^{2\alpha}} - \Gamma_1\ell^{-2\alpha}\|\theta\|_{L^{\infty}}^2\right)\le 0
\la{fina}
\ee
holds for any $t$, $x\in\Omega$ with $d(x)\ge \ell$ and $|h|\le\fr{\ell}{16}$.

We take $\delta>0$, $T>0$. We claim that, for any $\delta>0$ and any $T>0$
\[
\sup_{d(x)\ge \ell, |h|\le\fr{\ell}{16}, 0\le t \le T}\fr{|\delta_h\theta(x)|^2}{|h|^{2\alpha}} \le (1+\delta)\left[\|\theta_0\|_{C^{\alpha}}^2 + \Gamma_1 \ell^{-2\alpha}\|\theta_0\|_{L^{\infty}}^2\right] 
\]
holds.

The rest of the proof is done by contradiction. Indeed, assume by contradiction that there exists $\tilde{t}\le T$, $\tilde{x}$ and $\tilde{h}$ with $d({\tilde{x}})\ge \ell$ and $|{\tilde{h}}|\le {\fr{\ell}{16}}$ such that
\[
\fr{|\theta(\tilde{x} +\tilde{h}, \tilde{t})- \theta(\tilde{x},\tilde{t}) |^2}{|h|^{2\alpha}} >(1+\delta)\left[\|\theta_0\|_{C^{\alpha}}^2 + \Gamma_1 \ell^{-2\alpha}\|\theta_0\|_{L^{\infty}}^2\right] = R
\]
holds. Because the solution is smooth, we have
\[
\fr{|\delta_h\theta(x,t)|^2}{|h|^{2\alpha}} \le (1+\delta)\|\theta_0\|_{C^{\alpha}}^2
\]
for a short time $0\le t\le t_1$. (Note that this is not a statement about well-posedness in this norm: $t_1$ may depend on higher norms.) Also, because the solution is smooth, it is bounded in $C^1$, and
\[
\sup_{d(x)\ge \ell, |h|\le \fr{\ell}{16}}\fr{|\delta_h\theta(x)|^2}{|h|^2}\le C
\]
on the time interval $[0,T]$. It follows that there exists $\delta_1>0$ such that
\[
\sup_{d(x)\ge \ell, |h|\le \delta_1}\fr{|\delta_h\theta(x)|^2}{|h|^{2\alpha}}
\le C\delta_1^{2-2\alpha}\le \fr{R}{2}.
\]
In view of these considerations, we must have $\tilde{t} >t_1$, $|\tilde h|\ge\delta_1$. Moreover, the supremum is attained: there exists $\bar{x}\in\Omega$ with $d(\bar{x})\ge \ell$ and $\bar h \neq 0$ such that $\delta_1\le |\bar h|\le\fr{\ell}{16}$ such that
\[
\fr{|\theta(\bar{x}+\bar{h}, \tilde{t})-\theta(\bar{x},\tilde{t})|^2}{|\bar{h}|^{2\alpha}} =
s(\tilde{t}) =\sup_{d(x)\ge \ell, |h|\le \fr{\ell}{16}}\fr{|\delta_h\theta (\tilde{t})|^2}{|h|^{2\alpha}} > R.
\]
Because of (\ref{fina}) we have that
\[
\fr{d}{dt}\fr{|\theta(\bar{x}+\bar{h},t)-\theta(\bar{x},t)|^2}{|\bar{h}|^{2\alpha}}_{\left|\right. t=\tilde{t}} <0
\]
and therefore there exists $t'<\tilde{t}$ such that $s(t')>s(\tilde{t})$. This implies that
$\inf\{t>t_1\left|\right. s(t)>R\} = t_1$ which is absurd because we made sure that $s(t_1)< R$. Now $\delta$ and $T$ are arbitrary, so we have proved
\be
\sup_{d(x)\ge \ell, |h|\le\fr{\ell}{16}, t\ge 0}\fr{|\delta_h\theta(x)|^2}{|h|^{2\alpha}} \le \left[\|\theta_0\|_{C^{\alpha}}^2 + \Gamma_1 \ell^{-2\alpha}\|\theta_0\|_{L^{\infty}}^2\right]
\la{final}
\ee
which finishes the proof of the theorem.

\noindent{\bf Proof of Theorem \ref{alphaint}}. The proof follows from (\ref{final}) because $\Gamma_1$ does not depend on $\ell$. For any fixed $x\in\Omega$ we may take $\ell$ such that $\ell\le d(x)\le 2 \ell$. Then (\ref{final}) implies 
\be
d(x)^{2\alpha}\fr{|\delta_h\theta(x,t)|^2}{|h|^{2\alpha}}\le \left[\|\theta_0\|_{C^{\alpha}}^2 + \Gamma_1 2^{2\alpha}\|\theta_0\|_{L^{\infty}}^2\right].
\la{alpahintb}
\ee

\section{Gradient bounds}
We take the gradient of (\ref{sqg}). We obtain
\[
(\pa_t + u\cdot\na)\na\theta + (\na u)^*\na\theta +\na\l\theta =0
\] 
where $(\na u)^*$ is the transposed matrix. Let us take a good cutoff $\chi$. Then $g=\na\theta$ obeys everywhere
\be
\pa_t g + u\cdot\na g +\l(\chi g) + C_{\chi}(\theta) + (\na u)^*g = 0
\la{nateq}
\ee
with $C_{\chi}$ given in (\ref{cchi}). We multiply by $g$ and, using the fact that $\chi(x)=1 $ when $d(x)\ge \ell$ we obtain
\be
\fr{1}{2}L_{\chi}g^2 + D(f) + gC_{\chi}(\theta) + g(\na u)^*g = 0
\la{lchieqn}
\ee
when $d(x)\ge \ell$, where $L_{\chi}$ is similar to the one defined in (\ref{lchi}):
\be
L_{\chi}(\phi) = \pa_t\phi + u\cdot\na \phi + \l(\chi^2\phi)
\la{lchin}
\ee 
and $f=\chi g$.  Recall that $D(f) = f\l f-\l\left(\fr{f^2}{2}\right)$. Then, using (\ref{nachi}) and (\ref{nauxb}) we deduce
\be
\fr{1}{2}L_{\chi}g^2 + D(f) \le \fr{\Gamma_3}{d(x)^2}|g|\|\theta\|_{L^{\infty}(\Omega)}
+  C\left(\sqrt{\rho D(f)} + \|\theta\|_{L^{\infty}(\Omega)}\left(\fr{1}{d(x)} + \fr{1}{\rho}\right) + |\na\theta(x)|\right)g^2
\la{ginone}
\ee
for $d(x)\ge \ell$. Using a Young inequality we deduce
\be
L_{\chi}g^2 + D(f) \le  \fr{2\Gamma_3}{d(x)^2}\|\theta\|_{L^{\infty}(\Omega)}|g| 
+ C_4\rho g^4 + C_4\|\theta\|_{L^{\infty}(\Omega)}\left(\fr{1}{d(x)} + \fr{1}{\rho}\right)g^2 + C_4|g|^3
\la{gintwo}
\ee
for $d(x)\ge \ell$.  Now $|g| = |f|$ when $d(x)\ge \ell$. If $|g(x)|\ge M\|\theta\|_{L^{\infty}(\Omega)}d(x)^{-1}$ then, in view of (\ref{nlbd})
\be
D(f)\ge \gamma_2\|\theta\|_{C^{\alpha}(\Omega)}^{-\fr{1}{1-\alpha}}|g|^{3+\fr{\alpha}{1-\alpha}}(d(x))^{\fr{\alpha}{1-\alpha}} + \fr{\gamma_1}{d(x)}g^2
\la{dlow}
\ee
which is a super-cubic lower bound. We choose in this case
\be
\rho^{-1} = C_5 |g(x)|,
\la{rhoch}
\ee
and the right hand side of (\ref{gintwo}) becomes at most cubic in $g$: 
\be
L_{\chi}g^2 + D(f) \le  |g|^3\left[\fr{2\Gamma_3}{M^2\|\theta\|_{L^{\infty}(\Omega)}} + C_4\left (\fr{1}{C_5} + \fr{1}{M} + C_5\|\theta\|_{L^{\infty}(\Omega)} + 1\right)\right] = K|g|^3.
\la{ginthree}
\ee
In view of (\ref{dlow}) we see that
\be
L_{\chi}g^2 + |g|^3\left(\gamma_2\left(\|\theta\|_{C^{\alpha}(\Omega)}^{-\fr{1}{\alpha}}|g(x)|d(x)\right)^{\fr{\alpha}{1-\alpha}} - K\right)\le 0
\la{ginfour}
\ee
holds for $d(x)\ge \ell$, if $|g|\ge M\|\theta\|_{L^{\infty}}d(x)^{-1}$. 
In the opposite case, $|g(x)|\le M\|\theta\|_{L^{\infty}}d(x)^{-1}$ we choose
\be
\rho(x) = d(x)
\la{rhochoib}
\ee
and obtain from (\ref{gintwo})
\be
\ba
L_{\chi}g^2 + D(f) \\
\le \fr{1}{d(x)^3}\left[C_4M^4\|\theta\|_{L^{\infty}(\Omega)}^4 +   C_4M^3\|\theta\|_{L^{\infty}(\Omega)}^3 + 2C_4M^2\|\theta\|_{L^{\infty}(\Omega)}^3+ 2M\Gamma_3\|\theta\|_{L^{\infty}(\Omega)}^2\right] = \fr{K_1}{d(x)^3}
\ea
\la{gina}
\ee
and using the convex damping inequality (\ref{nlbd})
\[ 
D(f)\ge \gamma_1\fr{g^2}{d(x)}
\]
we obtain in this case
\be
L_{\chi}g^2 + \fr{1}{d(x)}\left(\gamma_1g^2(x) - \fr{K_1}{d(x)^2}\right)\le 0.
\la{gino}
\ee
Putting together (\ref{ginfour}) and (\ref{gino}) and \ref{hldrbnd} 
we obtain
\beg{thm}\la{gradb} Let $\theta$ be a smooth solution of (\ref{sqg}). Then
\be
\sup_{d(x)\ge \ell}|\na \theta(x,t)| \le C\left[\|\na\theta_0\|_{L^{\infty}(\Omega)} + \fr{P(\|\theta\|_{L^{\infty}(\Omega)})}{\ell}\right]
\la{fing}
\ee
where $P(\|\theta\|_{L^{\infty}(\Omega)})$ is a polynome of degree four.
\end{thm}

\noindent{\bf Proof of Theorem \ref{gradint}}. The proof follows by choosing $\ell$ depending on $x$, because the constants in (\ref{fing}) do not depend on $\ell$.
\section{Example: Half Space}
The  case of the half space is interesting because global smooth solutions of (\ref{sqg}) are easily obtained by reflection: If the initial data $\theta_0$ is smooth and compactly supported in  $\Omega = \Rr^d_+$ and if we consider its odd reflection
\be
\widetilde{\theta_0}(x) = \left\{
\begin{array}{c}
\theta_0(x_1,\dots x_d),\quad\quad \; {\mbox{if}}\; x_d>0,\\
-\theta_0(x_1,\dots, -x_d)\quad {\mbox{if}}\; x_d<0
\ea
\right.
\la{tildet}
\ee
then the solution of the critical SQG equation in the whole space, with intitial data $\widetilde{\theta_0}$ is globally smooth and its restriction to $\Omega$ solves (\ref{sqg}) there. This follows because of reflection properties of the heat kernel and of the Dirichlet Laplacian. 

The heat kernel with Dirichlet boundary conditions in $\Omega =\Rr^d_{+}$ is 
\[
H(x,y,t) = ct^{-\fr{d}{2}}\left( e^{-\fr{|x-y|^2}{4t}} - e^{-\fr{|x-{\widetilde {y}}|^2}{4t}}\right) 
\]
where $\widetilde{y} = (y_1,\dots, y_{d-1}, -y_d)$. More precisely,
\be
H(x,y,t) = G^{(d-1)}_t(x'-y')\left[G_t(x_d-y_d)-G_t(x_d+y_d)\right]
\la{Hsp}
\ee
with $x'= (x_1,\dots, x_{d-1})$,
\be
G_t^{(d-1)}(x') = \left(\fr{1}{4\pi t}\right)^{\fr{d-1}{2}}e^{-\fr{|x'|^2}{4t}}
\la{gprime}
\ee
and
\be
G_t(\xi) = \left(\fr{1}{4\pi t}\right)^{\fr{1}{2}} e^{-\fr{\xi^2}{4t}}
\la{G}
\ee
Let us note that
\be
\na_x H = H\left (
\begin{array}{c}
-\fr{x'-y'}{2t}\\
-\fr{x_d -y_d}{2t} + \fr{y_d}{t}\fr{e^{-\fr{x_dy_d}{t}}}{1-e^{-\fr{x_dy_d}{t}}}
\ea
\right)
\la{naxhh}
\ee
We check that (\ref{grbx}) is obeyed. Indeed, because $1-e^{-p}\ge \fr{p}{2}$ when $0\le p\le 1$ it follows that
\[
\fr{y_d}{t}e^{-\fr{x_dy_d}{t}}(1-e^{-\fr{x_dy_d}{t}})^{-1}\le \fr{y_d}{t}\fr{2t}{x_dy_d}
\]
if $\fr{x_dy_d}{t}\le 1$, and if $p=\fr{x_dy_d}{t} \ge 1$ then
\[
\fr{y_d}{t}e^{-\fr{x_dy_d}{t}}(1-e^{-\fr{x_dy_d}{t}})^{-1}\le \fr{e}{e-1}\fr{y_d}{t}e^{-\fr{x_dy_d}{t}}.
\]
In this case, if $\fr{x_d}{\sqrt{t}} \ge 1$ then $\fr{y_d}{t}\le t^{-\fr{1}{2}}p$ and $pe^{-p}$ is bounded; if $\fr{x_d}{\sqrt{t}}\le 1$ we write
$\fr{y_d}{t} = t^{-\fr{1}{2}}(\fr{y_d-x_d}{\sqrt{t}} + \fr{x_d}{\sqrt{t}})$ and thus we obtain:
\be
\left|\na_x H\right| \le C H\left [\fr{1}{\sqrt{t}}(1+ \fr{|x-y|}{\sqrt{t}}) + \fr{1}{x_d}\right]
\la{naxhhb}
\ee
We check (\ref{cancel1}): First we have
\be
(\na_x + \na_y)H = \left(
\begin{array}{c}
0\\
\fr{x_d+y_d}{t}G_t(x_d+y_d)G^{(d-1)}_t(x'-y')
\ea
\right)
\la{naxplusnay}
\ee
and then
\be
\int_{\Omega}\left| (\na_x+ \na_y)H(x,y,t)\right| dy \le Ct^{-\fr{1}{2}}e^{-\fr{x_d^{2}}{4t}}.
\la{hone}
\ee
Indeed, the only nonzero component occurs when the differentiation is with respect to the normal direction, and then
\be
\left|(\pa_{x_d} + \pa_{y_d})H(x,y,t)\right| = ct^{-\fr{d}{2}}e^{-\fr{|x'-y'|^2}{4t}}\left(\fr{x_d+y_d}{t}\right) e^{-\fr{(x_d+y_d)^2}{4t}}
\la{naxplus}
\ee
Therefore
\be
\ba
\int_{\Omega} \left|(\na_x+ \na_y)H(x,y,t)\right| dy \le Ct^{-\fr{1}{2}}\int_0^{\infty}\left(\fr{x_d+y_d}{t}\right )e^{-\fr{(x_d+y_d)^2}{4t}}dy_d\\
= Ct^{-\fr{1}{2}}\int_{\fr{x_d}{\sqrt {t}}}^{\infty}\xi e^{-\fr{\xi^2}{4}}d\xi
\\ 
= Ct^{-\fr{1}{2}}e^{-\fr{x_d^2}{4t}}.
\ea
\la{naxnaxc}
\ee
We check (\ref{cancel2}): first
\be
\ba
\pa_{x'}(\na_x  + \na_y)H = -\fr{x_d+y_d}{t}G_t(x_d+y_d)\fr{(x'-y')}{2t}G_t^{(d-1)}(x'-y')\\
\pa_{x_d}(\na_x + \na_y)H = \left(\fr{1}{t} - \fr{(x_d+y_d)^2}{2t^2}\right)G_t(x_d+y_d)G_t^{(d-1)}(x'-y')
\ea
\ee
Consequently
\be
|\na_x(\na_x + \na_y)H(x,y,t)|\le Ct^{-\fr{d}{2} -1}\left(1+\fr{|x'-y'|}{\sqrt{t}}\right)\left(1+ \fr{(x_d+y_d)^2}{t}\right)e^{-\fr{|x'-y'|^2}{4t}}e^{-\fr{(x_d+y_d)^2}{4t}}
\la{naxnaxplusnay}
\ee
and (\ref{cancel2}) follows:
\[
\int_{\Omega}|\na_x(\na_x + \na_y)H(x,y,t)|dy \le Ct^{-1}\int_{\fr{x_d}{\sqrt{2t}}}^{\infty}(1+z^2)e^{-\fr{z^2}{2}}dz.
\]

We compute $\Theta$ and $\l 1$:
\be
\Theta(x,t)= (e^{t\D}1)(x) = \int_{\Omega}H(x,y,t)dy = \fr{1}{\sqrt{2\pi}}\int_{-\fr{x_d}{\sqrt{2t}}}^{\fr{x_d}{\sqrt{2t}}}e^{-\fr{\xi^2}{2}}d\xi
\la{heat1}
\ee
and therefore 
\[
 \int_0^{\infty}t^{-\fr{3}{2}}(1-e^{t\D}1)dt =
\fr{2}{\sqrt{2\pi}}\int_0^{\infty}t^{-\fr{3}{2}}dt\int_{\fr{x_d}{\sqrt{2t}}}e^{-\fr{\xi^2}{2}}d\xi= \fr{4}{x_d\sqrt{\pi}}.
\]

\beg{rem}{\la{weak}} We note here that $\l^{s} 1 = C_sy_d^{-s}$ is calculated by duality:
\[
\ba
\left(\l^{s}1, \phi\right) = \left(1, \l^{s}\phi\right ) \\
=c_{s}\int_{\Omega}dx\int_0^{\infty}t^{-1-\fr{s}{2}}dt\left [\phi(x)-\int_{\Omega}H(x,y,t)\phi(y)dy\right]\\
=c_{s}\int_0^{\infty}t^{-1-\fr{s}{2}}dt\left[\int_{\Omega}\phi(x)dx -\int_{\Omega}\Theta(y_d,t)\phi(y)dy\right]\\
= c_{s}\int_0^{\infty}t^{-1-\fr{s}{2}}dt\int_{\Omega}\left(1-\Theta(y_d,t)\right)\phi(y)dy\\
=\fr{2c_{s}}{\sqrt{2\pi}}\int_{\Omega}\phi(y)\int_0^{\infty}t^{-1-\fr{s}{2}}dt\int_{\fr{y_d}{\sqrt{2t}}}^{\infty}e^{-\fr{\xi^2}{2}}d\xi\\
= C_{s}\int_{\Omega}y_d^{-s}\phi(y)dy
\ea
\]
where we used the symmetry of the kernel $H$ and (\ref{heat1}).

\end{rem}

We observe that if we consider horizontal finite differences, i.e. 
$h_d= 0$ then $C_h(\theta)$ vanishes,  and we deduce that 
\be
\sup_{x,h',t}|h'|^{-\alpha}|\theta (x'+ h', x_d, t)-\theta (x', x_d,t)|\le C_{1,\alpha} 
\la{partialh}
\ee
with $C_{1,\alpha}$ the partial $C^{\alpha}$ norm of the initial data.
This inequality can be used to prove that $u_2$ is bounded when $d=2$. Indeed
\be
u_2(x,t) = c\int_{\Omega}\left(\fr{1}{|x-y|^3} -\fr{1}{|x-\widetilde{y}|^3}\right)(x_1-y_1)\theta(y,t)dy
\la{uone}
\ee
and the bound is obtained using the partial H\"{o}lder bound on $\theta$ (\ref{partialh}) and the uniform bounds $\|\theta\|_{L^p}$ for $p=1, \infty$.
The outline of the proof is as follows: we split the integral
\be
u_2 = u_{2}^{in} + u_2^{out}
\la{uinout}
\ee
with
\be
u_2^{in}(x) = c\int_{|x_1-y_1|\le \delta, |x_2-y_2|\le \delta}\left(\fr{1}{|x-y|^3} -\fr{1}{|x-\widetilde{y}|^3}\right)(x_1-y_1)\left(\theta(y_1, y_2, t) - \theta( x_1, y_2,t)\right)dy
\la{uin}
\ee
and
\be
u_2^{out}(x) = c\int_{\max\{|x_1-y_1|, |x_2-y_2|\}\ge \delta}\left(\fr{1}{|x-y|^3} -\fr{1}{|x-\widetilde{y}|^3}\right)(x_1-y_1)\theta(y_1,y_2,t)dy
\la{uout}
\ee
where in (\ref{uin}) we used the fact that the kernel is odd in the first variable. Then, for $u^{in}$ we use the bound (\ref{partialh}) to derive
\be
|u_2^{in}(x)|\le C_{1,\alpha}C\int_0^{\sqrt{2}\delta}\rho^{-1+\alpha}d\rho =
CC_{1,\alpha}\delta^{\alpha}
\la{uinbound}
\ee
and for $u^{out}$,  if we have no other information on $\theta$ we just bound

\be
|u_2^{out}(x)|\le C\log\left(\fr{L}{\delta}\right )\|\theta_0\|_{L^{\infty}} + CL^{-2}\|\theta_0\|_{L^1}
\la{uoutbou}
 \ee
with some $L\ge \delta$. Both $\delta$ and $L$ are arbitrary. 

Finally, let us note that even if $\theta\in C_0^{\infty}(\Omega)$, the tangential component of the velocity need not vanish at the boundary because it is given by the integral
\[
u_1(x_1, 0, t) = -c\int_{\Rr^2_{+}}\fr{2y_2}{\left((x_1-y_1)^2 +y_2^2\right )^{\fr{3}{2}}} \theta(y,t)dy.
\]

\section{Appendix 1}
We sketch here the proofs of (\ref{naxnaxb}) (\ref{cancel1}) and (\ref{cancel2}). We take a point $\bar{x}\in\Omega$, a point $y\in\Omega$ and 
distinguish between two cases, if $d(\bar{x})< \fr{|\bar{x}-y|}{4}$ and if
$d(\bar{x})\ge \fr{|\bar{x}-y|}{4}$. In the first case we take a ball $B$ of radius $\delta=\fr{d(\bar{x})}{8}$ centered at $\bar{x}$ and in the second case we take also a ball $B$ centered at $\bar{x}$ but with radius $\delta=\fr{d(\bar{x})}{2}$. We note that in both cases the radius $\delta$ is proportional to $d(\bar{x})$. We take $x\in B(\bar{x}, \fr{\delta}{2})$, we fix $y\in\Omega$, take the function $h(z,t) = H_D(z,y,t)$, and apply Green's identity in the domain 
$U = B\times (0,t)$. We obtain 
\[
\ba
0=\int_U\left[(\pa_s-\Delta_z)h(z,s)G_{t-s}(x-z) + h(z,s)(\pa_s+\Delta_z)G_{t-s}(x-z)\right]dzds \\
= h(x,t)-G_t(x-y) + \int_0^t\int_{\pa B}\left[\fr{\pa G_{t-s}(x-z)}{\pa n}h(z,s)-\fr{\pa h(z,s)}{\pa n}G_{t-s}(x-z)\right]
\ea
\]
and thus
\[
H_D(x,y,t) = G_t(x-y) - \int_0^t\int_{\pa B}\left[\fr{\pa G_{t-s}(x-z)}{\pa n}h(z,s)-\fr{\pa h(z,s)}{\pa n}G_{t-s}(x-z)\right]
\]

We note that the $x$ dependence is only via $G$, and $x-z$ is bounded away from zero. 
We differentiate twice under the integral sign, and use the upper bounds (\ref{hb}), (\ref{grbx}). We have 
\[
\ba
|\na_x\na_x H_D(x,y,t)-\na_x\na_x G_{t}(x-y)| \\
\le C\int_0^t\int_{\pa B}(t-s)^{-\fr{d+3}{2}}p_3(\fr{|x-z|}{\sqrt{t-s}})e^{-\fr{|x-z|^2}{4(t-s)}}s^{-\fr{d}{2}}e^{-\fr{|y-z|^2}{Ks}}dzds\\
 +\int_0^{\min\{t; d^2(y)\}}\int_{\pa B} (t-s)^{-\fr{d+2}{2}}p_2(\fr{|x-z|}{\sqrt{t-s}})e^{-\fr{|x-y|^2}{4(t-s)}}s^{-\fr{d+1}{2}}p_1(\fr{|y-z|}{\sqrt{s}})e^{-\fr{|y-z|^2}{Ks}}dzds \\
+\int_{\min\{t; d^2(y)\}}^t\int_{\pa B} (t-s)^{-\fr{d+2}{2}}p_2(\fr{|x-z|}{\sqrt{t-s}})e^{-\fr{|x-y|^2}{4(t-s)}}s^{-\fr{d}{2}}\fr{1}{d(y)}\fr{w_1(y)}{|y-z|}e^{-\fr{|y-z|^2}{Ks}}dzds \\
\ea
\]
where $p_k(\xi)$ are polynomials of degree $k$.
The integrals are not singular. In both cases $|x-z|\ge \fr{\delta}{2}$,
and any negative power $(t-s)^{-\fr{k}{2}}$ can be absorbed by $e^{-\fr{|x-z|^2}{8(t-s)}}$ at the price $|x-z|^{-k}\le C\delta^{-k}$, still leaving  $e^{-\fr{|x-z|^2}{8(t-s)}}$ available. Similarly,
in the first case  $|y-z|\ge |\bar{x}-y| -\delta \ge \delta $ and in the second case $|y-z|\ge |\bar{x}-z| - |\bar{x}-y| \ge \fr{\delta}{2}$. Any power $s^{-\fr{k}{2}}$ can be absorbed by $e^{-\fr{|y-z|^2}{2Ks}}$ at the price $|y-z|^{-k}\le C\delta^{-k}$ still leaving $e^{-\fr{|y-z|^2}{2Ks}}$ available. We note that if $d(y)<d(x)$  so that $d(y)^2<t$ is possible, then, in view of (\ref{phione}) we have $\fr{w_1(y)}{|y-z|d(y)}\le C\delta^{-1}$.
 We also note that view of the fact that
\[
|x-y|^2t^{-1}\le 2\left(\fr{t-s}{t}\left(\fr{|x-z|^2}{t-s}\right) + \fr{s}{t}\left(\fr{|y-z|^2}{s}\right)\right)
\] 
we have a bound 
\[
e^{-\fr{|x-z|^2}{8(t-s)}-\fr{|y-z|^2}{2Ks}}\le e^{-\fr{|x-y|^2}{\tilde{K}t}}
\]
with $\tilde{K} = 16 + 4K$. Pulling this exponential out and estimating all the rest in terms of $\delta$ we obtain, in both cases, all the integrals bounded by
$Ct\delta^{-d-4}$ and therefore we have, in both cases,
\[
|\na_x\na_x H_D(x,y,t)-\na_x\na_x G_{t}(x-y)|\le Ce^{-\fr{|x-y|^2}{\tilde{K}t}}t\delta^{-d-4}\le Ct^{-1-\fr{d}{2}}e^{-\fr{|x-y|^2}{\tilde{K}t}}
\]
because $t\le c\delta^2$. This proves (\ref{naxnaxb}).

For (\ref{cancel1}) and (\ref{cancel2}) we start by noticing that it is enough to prove the estimates
\be
\int_{B(x, \fr{d(x)}{14})}|(\na_x+\na_y)H_D(x,y,t)|dy \le C t^{-\fr{1}{2}}e^{-\fr{d(x)^2}{Kt}}
\la{cone}
\ee
and
\be
\int_{B(x, \fr{d(x)}{14})}|\na_x(\na_x+\na_y)H_D(x,y,t)|dy \le C t^{-1}e^{-\fr{d(x)^2}{Kt}}
\la{ctwo}
\ee
for $t<cd^2(x)$. Indeed, if $|x-y|\ge \fr{d(x)}{14}$, individual Gaussian upper bounds for up to two derivatives of $H_D$ suffice (there is no need for cancellations).
In order to prove (\ref{cone}) and (\ref{ctwo}) we use a good cutoff $\chi$ with a scale $\ell = \fr{d(x)}{100}$. We take $y\in B(x, \fr{d(x)}{14})$. Both $x$ and $y$ are fixed for now. We note that the function
\[
z\mapsto h(z) = \chi(z)G_t(z-y)
\]
solves
\[
(\pa_t-\Delta)h(z,t) = -\left[(\Delta\chi(z))G_t(z-y) +2(\na\chi(z))\cdot\na G_t(z-y)\right] = F(z,y,t),
\]
vanishes for $z\in \pa\Omega$, and has initial datum $h_0= \chi(z)\delta(z-y)$, so, by Duhamel
\[
h(z,t) = e^{t\D}h_0 + \int_0^t e^{(t-s)\D}F(s)ds, 
\]
which, in view of $(e^{t\D}f)(z) = \int_{\Omega}H_D(z,w,t)f(w)dw$ yields
\[
\chi(z)G_t(z-y) = \chi(y)H_D(z,y,t) + \int_0^t\int_{\Omega}H_D(z,w,t-s)F(w,s)dwds
\]
for all $z$, and recalling that $\chi(x)=\chi(y) =1$, and reading at $z=x$ we have
\be
H_D(x,y,t) = G_t(x-y) + \int_0^t\int_{\Omega}H_D(x,w,t-s)\left[\Delta\chi(w)G_s(w-y) +2\na\chi(w)\cdot\na G_s(w-y)\right]dwds.
\la{repc}
\ee
The right hand side integral is not singular and can be differentiated  because the support of $\na\chi$ is far from the ball $B(x,\fr{d(x)}{14})$. Differentiation $\na_x+\na_y$ cancels the Gaussian $G_t(x-y)$. The estimates of the right hand side
\[
\left|(\na_x+\na_y) \int_0^t \int_{\Omega}H_D(x,w,t-s)F(w,y,s)dwds\right| \le Ct^{-\fr{d+1}{2}}e^{-\fr{d(x)^2}{Kt}}
\]
and
\[
\left|\na_x(\na_x+\na_y) \int_0^t \int_{\Omega}H_D(x,w,t-s)F(w,y,s)dwds\right| \le Ct^{-{\fr{d+2}{2}}}e^{-\fr{d(x)^2}{Kt}}
\]
for $t<cd^2(x)$ follow from Gaussian upper bounds. Integration $dy$ on the ball
$B(\fr{d(x)}{14})$ picks up the volume of the ball, and thus (\ref{cone}) and
(\ref{ctwo}) are verified.

\section{Appendix 2}
We sketch here the proof of local wellposedness of the equation (\ref{sqg}).
We start by defining a Galerkin approximation. We consider the projectors
$P_n$
\be
P_n f= \sum_{j=1}^nf_jw_j
\la{pn}
\ee
with $f_j= \int_{\Omega}f(x)w_j(x)dx$. We consider for fixed $n$ the approximate system
\be
\pa_t \theta_n + P_n\left(u_n\cdot\na\theta_n\right) + \l\theta_n = 0
\la{sqgn}
\ee
where
\be
u_n = \na^{\perp}\l^{-1}\theta_n = R_D^{\perp}\theta_n
\la{un}
\ee
with 
\be
(P_n\theta_n)(x,t) = \theta_n(x,t) = \sum_{j=1}^n\theta_{n,j}(t)w_j(x)
\la{pnt}
\ee
and with initial data $\theta_n(0) = P_n\theta_0$ where $\theta_0$ is a fixed smooth function belonging to $H_0^1(\Omega)\cap H^2(\Omega)$.  Although it was written as a PDE, the system (\ref{sqgn}) is a system of ODEs for the coefficients $\theta_{n,j}(t)= \int_{\Omega}\theta_nw_jdx$. Let us note that
$P_n$ does not commute with $\na$ but does commute with $-\Delta$ and functions of it. The function $u_n$ is divergence-free and it is a finite sum of divergence-free functions, 
\be
u_n(x) = \sum_{j=1}^n\lambda_j^{-\fr{1}{2}}\theta_{n,j}(t)\na^{\perp}w_j(x). 
\la{unwj}
\ee
Note however that $u_n\notin P_nL^{2}(\Omega)$. The normal component of $u_n$ vanishes at the boundary because $\na^{\perp}w_j\cdot\nu_{\left|\right.\Omega} =0$.
Moreover, because 
\[
\int_{\Omega}P_n(u_n\cdot\na\theta_n)\theta_n dx = \int_{\Omega}(u_n\cdot\na\theta_n)\theta_ndx = 0
\]
it follows that $\|\theta_n(t)\|_{L^2(\Omega)}$ is bounded in time and therefore the solution exists for all time.  The following upper bound for higher norms is uniform only for short time, and it is the bound that is used for local existence of smooth solutions. We apply $\l^2=-\Delta$ to (\ref{sqgn}) and use the fact that it is a local operator, it commutes with $P_n$ and with derivatives:
\be
\pa_t\l^2\theta_n + P_n\left(u_n\cdot\na\l^2\theta_n -2\na u_n\na\na\theta_n + (\l^2 u_n)\cdot\na\theta_n\right) + \l^3\theta_n = 0
\la{lapsqgn}
\ee
We take the scalar product with $\l^2\theta_n$. Because this is finite linear combinations of eigenfunctions, it vansihes at $\pa\Omega$ and integration by parts is allowed. We obtain
\be
\ba
\fr{d}{2dt}\|\l^2\theta_n\|^2_{L^2(\Omega)} + \|\l^{\fr{5}{2}}\theta_n\|^2_{L^2(\Omega)}\\
\le\|\l^2u_n\|^2_{L^2(\Omega)}\|\l^2\theta_n\|_{L^2(\Omega)}\|\na\theta_n\|_{L^{\infty}(\Omega)} + 2\|\na u_n\|_{L^{\infty}(\Omega)}\|\na\na\theta\|_{L^2(\Omega)}\|\l^2\theta\|_{L^2(\Omega)} 
\ea
\la{incomplete}
\ee
We note now that
\be
\l^2u_n = \sum_{j=1}^n\theta_{n,j}(-\D)\lambda_j^{-\fr{1}{2}}\na^{\perp}w_j=
\na^{\perp}\l^{-1}(\l^2\theta_n) = R_D^{\perp}(\l^2\theta_n).
\la{delun}
\ee
Now $R_D$ is bounded in $L^2(\Omega)$ (It is in fact an isometry on components; this follows from (\ref{kat})),  therefore 
\be
\|\l^2 u_n\|_{L^2{(\Omega)}} \le \|\l^2\theta_n\|_{L^2(\Omega)}.
\la{delunb}
\ee
The fact that $R_D$ is bounded in $L^4(\Omega)$ is also true (\cite{jerisonkenig}). Then
\be
\|\l^2 u_n\|_{L^4{(\Omega)}} \le \|\l^{2}\theta_n\|_{L^4(\Omega)}.
\la{dellunb}
\ee
Moreover, it is known (see for instance (\cite{cabre})) that in $d=2$ we have
\[
\|f\|_{L^4(\Omega)} \le C\|\l^{\fr{1}{2}}f\|_{L^2(\Omega)}
\]
and therefore
\be
\|\D \theta_n\|_{L^{4}(\Omega)}\le \|\l^{\fr{5}{2}} \theta_n\|_{L^2{(\Omega)}}.
\la{delfourthn}
\ee
and
\be
\|\D u_n\|_{L^{4}(\Omega)}\le C\|\l^{\fr{5}{2}} \theta_n\|_{L^2{(\Omega)}}.
\la{delfourun}
\ee
Now we use the Sobolev embedding
\be
\|\na \phi\|_{L^{\infty}(\Omega)}\le C\left(\|\D\phi\|_{L^{4}(\Omega)} + \|\phi\|_{L^2(\Omega)}\right)
\la{sob}
\ee
and deduce, using also a Poincar\'{e} inequality 
\be
\fr{d}{dt}\|\l^2\theta_n\|_{L^2(\Omega)}^2 + \|\l^{\fr{5}{2}}\theta_n\|^2_{L^2(\Omega)} \le C\|\l^2\theta_n\|_{L^2(\Omega)}^2\|\l^{\fr{5}{2}}\theta\|_{L^2(\Omega)}.
\la{complete}
\ee
Thus, after a Young inequality we deduce that
\be
\sup_{t\le T} \|\l^2\theta_n\|_{L^2(\Omega)}^2 + \int_0^T\|\l^{\fr{5}{2}}\theta_n\|^2_{L^2(\Omega)}dt \le C\|\l^2\theta_0\|_{L^2(\Omega)}^2
\la{finn}
\ee
holds for $T$ depending only on $\|\l^2\theta_0\|_{L^2(\Omega)}$, with a constant independent of $n$. The following result can now be obtained by assing to the limit in a subsequence and using a Aubin-Lions lemma (\cite{lions}): 
\beg{prop} Let $\theta_0\in H_0^1(\Omega)\cap H^{2}(\Omega)$ in $d=2$. There exists $T>0$ a unique solution of (\ref{sqg}) with initial datum $\theta_0$ satisfying
\be
\theta\in L^{\infty}(0,T; H_0^1(\Omega)\cap H^{2}(\Omega))\cap L^2\left(0,T; {\mathcal{D}}\left(\l^{2.5}\right) \right).
\la{loc}
\ee
\end{prop}
Higher regularity can be obtained as well. Because the proof uses $L^2$- based Sobolev spaces and Sobolev embedding, it is dimension dependent. A proof in higher dimensions is also possible, but it requires using higher powers of $\D$, and will not be pursued here.

{\bf{Acknowledgment.}} The work of PC was partially supported by  NSF grant DMS-1209394


\begin{thebibliography}{99}
\bibitem{caf} L.A. Caffarelli and A.~Vasseur,
\newblock Drift diffusion equations with fractional diffusion and the
  quasi-geostrophic equation,  Ann. of Math., 171(3) (2010), 1903--1930.

\bibitem{cabre} X. Cabre, J. Tan, Positive solutions of nonlinear problems involving the square root of the Laplacian, Adv. Math. 224 (2010), no. 5, 2052-2093.
\bibitem{castro} A. Castro, D. C\'{o}rdoba, J. Gomez-Serrano and A. Martin,
Remarks on geometric properties of the SQG sharp front and the alpha-patches, Discrete Continuos Dynamical Systems  34 (12), (2014) 5045-5059.

\bibitem{cascor} A. Castro, D. ~C\'{o}rdoba, J. Gomez-Serrano, Global smooth solutions for the inviscid SQG equation (with . Preprint arxiv:1603.03325.

\bibitem{c} P. Constantin, Geometric statistics in turbulence, SIAM Review 
{\bf{36}}, (1) (1994). 73-98. 
\bibitem{ccw}
P.~Constantin, D.~C{\'o}rdoba, and J.~Wu. On the critical dissipative quasi-geostrophic equation,  {Indiana Univ. Math. J.}, 50(Special Issue):Dedicated to Professors Ciprian Foias and Roger Temam (Bloomington,IN, 2000). (2001), 97-107.


\bibitem{ci} P. Constantin, M. Ignatova, Remarks on the fractional Laplacian with Dirichlet boundary conditions and applications, IMRN, 2016.

\bibitem{cnum}P.~Constantin, M.-C. Lai, R.~Sharma, Y.-H. Tseng, and J.~Wu, New numerical results for the surface quasi-geostrophic equation, J. Sci. Comput., 50(1) (2012), 1--28.


\bibitem{cmt} P.~Constantin, A.J. Majda, and E.~Tabak, Formation of strong fronts in the {$2$}-{D} quasigeostrophic thermal active scalar,  Nonlinearity, 7(6) (1994), 1495--1533.


\bibitem{cv1} P. Constantin, V. Vicol, Nonlinear maximum principles for dissipative linear nonlocal operators and applications, GAFA {\bf{22}} (2012) 1289-1321.

\bibitem{cvt} P. Constantin, A. Tarfulea, V. Vicol, Long time dynamics of forced critical SQG, Communications in Mathematical Physics 335 (2015), no. 1, 93-141. 


\bibitem {cw}
P.~Constantin and J.~Wu, Regularity of {H}{\"o}lder continuous solutions of the supercritical quasi-geostrophic equation, Ann. Inst. H. Poincar{\'e} Anal. Non Lin{\'e}aire, 25(6) (2008) 1103--1110.


\bibitem{cord}D.~C{\'o}rdoba, Nonexistence of simple hyperbolic blow-up for the quasi-geostrophic equation, Ann. of Math. (2), 148(3) (1998) 1135--1152.


\bibitem{cc} A.~C\'ordoba, D.~C\'ordoba, { A maximum principle applied to quasi-geostrophic equations}. Comm.~Math.~Phys.~{\bf 249} (2004), 511--528.

\bibitem{cornum} D.~C{{\'o}}rdoba, M.A. Fontelos, A.M. Mancho, and J.L. Rodrigo, Evidence of singularities for a family of contour dynamics equations, Proc. Natl. Acad. Sci. USA, 102(17) (2005), 5949--5952.


\bibitem{davies1} E.B. Davies, Explicit constants for Gaussian upper bounds on heat kernels, Am. J. Math {\bf{109}} (1987) 319-333.


\bibitem{fr} C.~Fefferman and J.L. Rodrigo.  Analytic sharp fronts for the surface quasi-geostrophic equation. Comm. Math. Phys., 303(1) (2011), 261--288.

\bibitem{held} I.M. Held, R.T. Pierrehumbert, S.T. Garner, and K.L. Swanson, Surface quasi-geostrophic dynamics, J. Fluid Mech., 282 (1995),1--20.


\bibitem{jerisonkenig} D. Jerison, C. Kenig, The inhomogeneous Dirichlet problem in Lipschitz domains, J. Funct. Analysis {\bf{130}} (1995), 161-212.

\bibitem{knv} A.~Kiselev, F.~Nazarov, and A.~Volberg, Global well-posedness for the critical 2{D} dissipative quasi-geostrophic equation, Invent. Math., 167(3) (2007), 445--453.

\bibitem{lions} J.L. Lions, Quelque methodes de r\'{e}solution des probl\`{e}mes aux limites non lin\'{e}aires, Paris, Dunod (1969). 

\bibitem{res}S.G. Resnick, Dynamical problems in non-linear advective partial differential equations, ProQuest LLC, Ann Arbor, MI, 1995, Thesis (Ph.D.)--The University of Chicago.



\bibitem{svz} L.~Silvestre, V.~Vicol, and A.~Zlato{\v{s}}, On the {L}oss of {C}ontinuity for {S}uper-{C}ritical {D}rift-{D}iffusion {E}quations, Arch. Ration. Mech. Anal., 207(3) (2013) 845--877.



\bibitem{qszhang1} Q. S. Zhang, The boundary behavior of heat kernels of Dirichlet Laplacians, J. Diff. Eqn {\bf{182}} (2002), 416-430.

\bibitem{qszhang2} Q. S. Zhang, Some gradient estimates for the heat equation on domains and for an equation by Perelman, IMRN (2006), article ID92314, 1-39.

\end{thebibliography}
\end{document}